\documentclass[final]{aims}
\usepackage{amsmath,amsthm,amssymb,amsfonts}
\usepackage{mathtools}
\usepackage[utf8]{inputenc} 
\usepackage[T1]{fontenc} 
\usepackage{enumerate}
\usepackage{pdfpages}
\usepackage{graphicx}
\usepackage{tabularx}
\usepackage{hyperref}
\usepackage{color}
\usepackage{url}
\hypersetup{urlcolor=black,citecolor=black}
\def\currentvolume{X}
\def\currentissue{X}
\def\currentyear{200X}
\def\currentmonth{XX}
\def\ppages{X--XX}
\def\DOI{10.3934/xx.xx.xx.xx}
\renewcommand{\paragraph}[1]{\emph{#1.}}
\newcommand{\NN}{\mathbb{N}}

\newcommand{\RR}{\mathbb{R}}
\newcommand{\CC}{\mathbb{C}}
\newcommand{\dd}{\text{d}}
\newcommand{\ee}{\text{e}}
\newcommand{\ii}{\text{i}}
\newcommand{\abs}[1]{|#1|}
\newcommand{\absbig}[1]{\big|#1\big|}

\newcommand{\nH}{\mathcal{H}}
\newcommand{\nJ}{\mathcal{J}}
\newcommand{\nL}{\mathcal{L}}
\newcommand{\nO}{\mathcal{O}}
\renewcommand{\labelenumii}{\labelenumi\arabic{enumii}. }
\title[Gross--Pitaevskii equations with rotation term]{Efficient time integration methods for \\ Gross--Pitaevskii equations with rotation term} 
\author[P.~Bader, S.~Blanes, F.~Casas, M.~Thalhammer]{}
\email{bader@uji.es}
\email{serblaza@imm.upv.es}
\email{fernando.casas@uji.es}
\email{mechthild.thalhammer@uibk.ac.at}
\thanks{Version of \today}
\begin{document}
\maketitle 
$\,$ \\[-8mm]
\centerline{\scshape Philipp Bader} \\
\centerline{\footnotesize Universitat Jaume~I, Departament de Matem{\`a}tiques, 12071~Castell\'on, Spain.} \\[2mm]
\centerline{\scshape Sergio Blanes} \\
\centerline{\footnotesize Universitat Polit{\`e}cnica de Val{\`e}ncia, Instituto Universitario de Matem{\'a}tica Multi\-disciplinar, 46022~Valencia, Spain.} \\[2mm]
\centerline{\scshape Fernando Casas} \\
\centerline{\footnotesize Universitat Jaume~I, IMAC and Departament de Matem{\`a}tiques, 12071~Castell\'on, Spain.} \\[2mm]
\centerline{\scshape Mechthild Thalhammer} \\
\centerline{\footnotesize Leopold--Franzens--Universit{\"a}t Innsbruck, Institut f{\"u}r Mathematik, 6020~Innsbruck, Austria.} 
\begin{abstract}
The objective of this work is the introduction and investigation of favourable time integration methods for the Gross--Pitaevskii equation with rotation term.
Employing a reformulation in rotating Lagrangian coordinates, the equation takes the form of a nonlinear Schr{\"o}dinger equation involving a space-time-dependent potential.  
A natural approach that combines commutator-free quasi-Magnus exponential integrators with operator splitting methods and Fourier spectral space discretisations is proposed.  
Furthermore, the special structure of the Hamilton operator permits the design of specifically tailored schemes. 
Numerical experiments confirm the good performance of the resulting exponential integrators.
\end{abstract}
\section{Introduction}
\label{sec:Section1}
$\,$ \\[-1mm]

\paragraph{Scope}
Our main concern is the introduction and investigation of efficient numerical methods for nonlinear evolution equations involving an explicit time-dependency
\begin{equation} 
\label{eq:IntroductionEqu1} 
u'(t) = F\big(t, u(t)\big)\,, \quad t \in (t_0, T)\,; 
\end{equation}
previous work on non-autonomous linear evolution equations 
\begin{equation}
\label{eq:IntroductionEqu2} 
u'(t) = F(t) \, u(t)\,, \quad t \in (t_0, T)\,,
\end{equation}
has enhanced our interest in this subject.
In the present work, we center our attention on time-dependent Gross--Pitaevskii equations with rotation term modelling the creation of vortices in Bose--Einstein condensates. 

\paragraph{Exponential integrators}
Theoretical and numerical studies of time integration methods for systems of non-autonomous linear ordinary differential equations that can be cast into the form~\eqref{eq:IntroductionEqu2} with matrix-valued time-dependent function $F: [t_0, T] \to \CC^{M \times M}$ are found in~\cite{AlvermannFehske2011,AlvermannFehskeLittlewood2012,BaderBlanesKopylov2018, BlanesCasasGonzalezThalhammer2018,BlanesCasasThalhammer2017,BlanesCasasThalhammer2018}, see also references given therein.
Numerical comparisons presented within the context of large-scale applications arising in quantum physics clearly indicate that exponential integrators based on the Magnus expansion are favourable compared to standard integrators such as explicit Runge--Kutta methods;
moreover, it is demonstrated that an approach which avoids the actual computation of matrix commutators
\begin{equation*}
\big[F(s), F(t)\big] = F(s) \, F(t) - F(t) \, F(s)\,, \quad s, t \in [t_0, T]\,,
\end{equation*}
by considering compositions of matrix exponentials leads to approximations that are superior to those obtained by Magnus integrators.
In view of an appropriate extension to partial differential equations and finally to the nonlinear case~\eqref{eq:IntroductionEqu1}, it is helpful to understand this approach as a replacement of the exact evolution operator associated with~\eqref{eq:IntroductionEqu2} by a sequence of autonomous evolution equations, which are essentially of the form 
\begin{equation}
\label{eq:IntroductionEqu3}
u'(t) = F(t_*) \, u(t)\,, \quad t \in (t_0, T)\,,
\end{equation}
for certain fixed times $t_{*} \in [t_0, T]$;
we refer to this type of time integration methods as commutator-free quasi-Magnus (CFQM) exponential integrators.

\paragraph{Computational aspects}
For non-autonomous linear ordinary differential equations or spatial semi-discretisations of non-autonomous linear partial differential equations leading to ordinary differential equations of the form~\eqref{eq:IntroductionEqu2}, respectively, the realisation of CFQM exponential integrators relies on the numerical computation of matrix exponentials.
For non-autonomous linear Schr{\"o}dinger equations of a special structure
\begin{subequations}
\label{eq:IntroductionEqu4}
\begin{equation}
\ii \, u'(t) = \ii \, F(t) \, u(t) = H(t) \, u(t) = \big(A + B(t)\big) \, u(t)\,, \quad t \in (t_0, T)\,, 
\end{equation}
written as abstract evolution equations, an alternative methodology is to employ operator splitting~\cite{BlanesMoan2002Splitting,McLachlanQuispel2002} combined with spectral methods;
for instance, in situations where the Hamiltonian comprises the Laplace operator with respect to the spatial variables and a real-valued space-time-dependent potential $W: \RR^d \times [t_0, T] \to \RR$ that satisfies appropriate regularity and periodicity conditions
\begin{equation}
H(t) = - \, \tfrac{1}{2} \, \Delta + W(\cdot, t)\,, \quad A = - \, \tfrac{1}{2} \, \Delta\,, \quad B(t) = W(\cdot, t)\,,
\end{equation}
\end{subequations}
favourable approximations to~\eqref{eq:IntroductionEqu3} are based on suitable compositions of the solutions to the associated autonomous subequations
\begin{equation*}
\ii \, u'(t) = A \, u(t)\,, \quad \ii \, u'(t) = B(t_*) \, u(t)\,,
\end{equation*}
computed by Fourier spectral space discretisation and pointwise multiplication, respectively.

\paragraph{Modified integrators}
For Schr{\"o}dinger equations defined by a Hamilton operator of the special structure~\eqref{eq:IntroductionEqu4}, the incorporation of matrix commutators corresponds to modifications of the potential; 
indeed, straightforward calculations show that certain Lie-brackets of the second space derivative and a multiplication operator are given by
\begin{equation}
\label{eq:IntroductionLieCommutator}
\begin{gathered}
\big[\partial_{xx}, W\big] w = \partial_{xx} \big(W w\big) - W \partial_{xx} w = \big(\partial_{xx} W + 2 \, \partial_{x} W \, \partial_{x}\big) \, w\,, \\
\big[[\partial_{xx}, W], W\big] w
= 2 \, \partial_{x} W \, \partial_{x} (W w) - 2 \, W \, \partial_{x} W \, \partial_{x} w
= 2 \, \big(\partial_{x} W\big)^2 w\,.
\end{gathered}
\end{equation}
Together with Taylor series expansions and polynomial interpolation, this permits to design higher-order schemes with improved efficiency, see~\cite{BaderBlanesKopylov2018}; 
we refer to this type of methods as modified CFQM exponential integrators, and with regard to the initials of the authors, to a specific sixth-order scheme as BBK scheme. 
In combination with operator splitting, this approach leads to autonomous subequations of the form 
\begin{equation*}
\ii \, u'(t) = A \, u(t)\,, \quad \ii \, u'(t) = \widehat{B} \, u(t)\,, \quad \ii \, u'(t) = \big(\widehat{B} + \widetilde{B}\big) \, u(t)\,,  
\end{equation*}
where~$\widehat{B}$ represents certain linear combinations of the values of the potential and~$\widetilde{B}$ involves first space derivatives of the potential;
again, their numerical solution relies on Fourier spectral space discretisation and pointwise multiplication, respectively.

\paragraph{Extension to nonlinear equations}
In this work, we provide an expedient extension of CFQM exponential integrators and modified CFQM exponential integrators for non-autonomous linear evolution equations~\eqref{eq:IntroductionEqu2} to the substantially more complex nonlinear case~\eqref{eq:IntroductionEqu1}.
With regard to relevant quantum physical applications in the field of Bose--Einstein condensation, modelled by time-dependent Gross--Pitaevskii equations with rotation term, we focus on non-autonomous nonlinear evolution equations of Schr{\"o}dinger type that can be cast into the form 
\begin{equation}
\label{eq:IntroductionEqu5}
\ii \, u'(t) = \ii \, F\big(t, u(t)\big) = H\big(t, u(t)\big) = \Big(A + B(t) + C\big(u(t)\big)\Big) \, u(t)\,, \quad t \in (t_0, T)\,; 
\end{equation}
comparably to~\eqref{eq:IntroductionEqu4}, $A$ represents the Laplace operator, $B$ a real-valued space-time-dependent trapping potential, and~$C$ a real-valued nonlinear multiplication operator.
The restriction to this particular setting permits the implementation of CFQM exponential integrators as well as modified CFQM exponential integrators
by operator splitting and Fourier spectral methods.

\paragraph{Strategy for full discretisation}
We point out that the extension of the afore mentioned concepts to the nonlinear case is not entirely straightforward, since different approaches seem likely and the overall computational cost has to be taken into account.
\begin{enumerate}[(i)]
\item
A strategy for the discretisation of~\eqref{eq:IntroductionEqu5}, which is in principle practicable, uses operator splitting, that is, a decomposition into linear and nonlinear terms 
\begin{equation*}
\ii \, u'(t) = A \, u(t) + B(t) \, u(t)\,, \quad \ii \, u'(t) = C\big(u(t)\big) \, u(t)\,.
\end{equation*}
The numerical solution of the non-autonomous linear subequation relies on CFQM exponential integrators and once again on operator splitting methods;
due to an invariance property that is characteristic for nonlinear Schr{\"o}dinger equations, the nonlinear subequation reduces to a linear equation
\begin{equation*}
\ii \, u'(t) = C\big(u(t_{*})\big) \, u(t)\,,
\end{equation*}
resolvable by pointwise multiplication. 
For a splitting method comprising~$r$ compositions and a CFQM exponential integrator involving~$s$ compositions, this strategy requires the numerical solution of $r^2 s$ linear subequations by Fourier spectral space discretisation as well as $r + r^2 s$ linear subequations by pointwise multiplication, since 
\begin{equation*}
\begin{gathered}
\text{application of $r$-stage splitting method:} \\
\begin{cases}
\text{$r$ times:} \quad &\ii \, u'(t) = A \, u(t) + B(t) \, u(t)\,, \\
\text{$r$ times (pointwise multiplication):} \quad &\ii \, u'(t) = C\big(u(t_{*})\big) \, u(t)\,, 
\end{cases} \\
\text{application of $s$-stage CFQM exponential integrator:} \\
\begin{cases}
\text{$r s$ times:} \quad \ii \, u'(t) = A \, u(t) + B(t_{*}) \, u(t)\,, 
\end{cases} \\
\text{application of $r$-stage splitting method:} \\
\begin{cases}
\text{$r^2 s$ times (Fourier spectral method):} \quad &\ii \, u'(t) = A \, u(t)\,, \\
\text{$r^2 s$ times (pointwise multiplication):} \quad &\ii \, u'(t) = B(t_{*}) \, u(t)\,.
\end{cases} 
\end{gathered}
\end{equation*}
\item
In view of CFQM exponential integrators and splitting methods involving a higher number of compositions, potentially yielding higher accuracy and efficiency, it is advantageous to employ a different strategy. 
At first, a CFQM exponential integrator is applied, which leads to a sequence of autonomous nonlinear evolutions equations
\begin{equation*}
\ii \, u'(t) = A \, u(t) + B(t_{*}) \, u(t) + C\big(u(t)\big) \, u(t)\,, 
\end{equation*}
in accordance with~\eqref{eq:IntroductionEqu3};
combined with operator splitting, this alternative approach requires the numerical solution of linear subequations of the form 
\begin{equation*}
\ii \, u'(t) = A \, u(t)\,, \quad \ii \, u'(t) = B(t_{*}) \, u(t) + C\big(u(s_{*})\big) \, u(t)\,,
\end{equation*}
for certain fixed values $s_{*}, t_{*} \in [t_0, T]$ by Fourier spectral space discretisation and pointwise multiplication, respectively, each~$r s$ times, since 
\begin{equation*}
\begin{gathered}
\text{application of $s$-stage CFQM exponential integrator:} \\
\begin{cases}
\text{$s$ times:} \quad \ii \, u'(t) = A \, u(t) + B(t_{*}) \, u(t) + C\big(u(t)\big) \, u(t)\,, 
\end{cases} \\
\text{application of $r$-stage splitting method:} \\
\begin{cases}
\text{$r s$ times (Fourier spectral method):} \quad &\ii \, u'(t) = A \, u(t)\,, \\
\text{$r s$ times (pointwise multiplication):} \quad &\ii \, u'(t) = B(t_{*}) \, u(t) + C\big(u(s_{*})\big) \, u(t)\,.
\end{cases} 
\end{gathered}
\end{equation*}
\end{enumerate}

\paragraph{Magnus integrators}
It is also noteworthy that a formal extension of standard Magnus integrators would lead to involved and impractical numerical approximations. 
More precisely, rewriting the non-autonomous nonlinear differential equation~\eqref{eq:IntroductionEqu1} by means of the formal calculus of Lie-derivatives as non-autonomous linear differential equation
\begin{equation*}
u'(t) = L(t) \, u(t)\,, \quad L(t) \equiv F\big(t, u(t)\big) \frac{\delta}{\delta u}\,, \quad t \in (t_0, T)\,,
\end{equation*}
the associated solution representation based on the Magnus expansion
\begin{equation*}
\begin{gathered}
u(t) = \ee^{\Omega(t - t_0)} \, u_0\,, \\
\Omega(t - t_0)
= \int_{t_0}^{t} L(\tau) \; \dd \tau + \tfrac{1}{2} \int_{t_0}^{t} \int_{t_0}^{\tau} \big[L(\tau), L(\sigma)\big] \; \dd \sigma \, \dd \tau + \dots, \\
t \in [t_0, T]\,,
\end{gathered}
\end{equation*}
involves nested integrals of Lie-brackets.
Specifically, for non-autonomous nonlinear Schr{\"o}dinger equations~\eqref{eq:IntroductionEqu5}, the Lie-bracket~$[L(\tau), L(\sigma)]$ corresponds to a second-order differential operator with coefficients depending on the solution values.

\paragraph{High accuracy}
A characteristic of rotating Bose--Einstein condensates, modelled by time-dependent Gross--Pitaevskii equations with rotation term, is the creation of vortices;
in order to simulate this quantum physical phenomenon correctly, high spatial and temporal resolution is required, and thus it is of relevance to employ accurate and efficient full discretisation methods.
Given the relevance of the theme, numerous contributions are concerned with the development of advanced methodologies and illustrate their benefits over standard numerical methods; 
we in particular mention an approach based on generalised-Laguerre–Fourier–Hermite spectral space discretisation and operator splitting methods that has been studied in~\cite{BaoLiShen2009, HofstaetterKochThalhammer2014} and the transformation to the rotating frame combined with the application of exponential integrators as proposed in~\cite{AntoineBesseBao2013,Bader2014,Bader2015,BaoMarahrensTangZhang2013,BesseDujardinLacroix2017}.
We point out that a rigorous stability and local error analysis has shown the reliability of generalised-Laguerre–Fourier–Hermite pseudo-spectral time-splitting methods, see~\cite{HofstaetterKochThalhammer2014};
but, due to the fact that the implementation of these specialised spectral methods is involved and computationally expensive, we believe that it is expedient to develop~\cite{BaoMarahrensTangZhang2013} further. 
In order to design higher-order schemes with improved accuracy and efficiency, we transform the considered time-dependent Gross--Pitaevskii equation to rotating Lagrangian coordinates;
this yields a non-autonomous nonlinear Schr{\"o}dinger equation of the form~\eqref{eq:IntroductionEqu5}, where the afore sketched strategy can be realised by fast Fourier transforms.

\paragraph{Notation}
The use of standard notation implicates that the same letter denotes diverse quantities in different contexts. 
For instance, in Sections~\ref{sec:Section2} and~\ref{sec:Section5}, the real number $\Omega \in \RR$ denotes the constant angular velocity;
by contrast, in Section~\eqref{sec:Section3}, the function $\Omega: [t_0, t_0 + h] \to \RR$ denotes the exponent in the Magnus expansion.
\section{Rotational Gross--Pitaevskii equation}
\label{sec:Section2}
$\,$ \\[-1mm]

\paragraph{Rotating Bose--Einstein condensates}
The phenomenon of Bose--Einstein condensation, for the first time observed in~1995~\cite{Anderson1995,Bradley1995,Davis1995}, has opened new perspectives on quantum physics. 
Over the last two decades, rotating Bose--Einstein condensates with the creation of quantised vortices as a particular indication of superfluidity~\cite{Abo2001,Madison2001} have attracted a lot of interest; 
various theoretical and experimental studies make a contribution to a deeper understanding of the observed quantum physical effects. 

\paragraph{Rotational Gross--Pitaevskii equation}
At temperatures significantly below the critical temperature, a prevalent model for the dynamical behaviour of a rotating Bose--Einstein condensate is the time-dependent Gross--Pitaevskii equation with additional angular momentum rotation term.
We study the following normalised formulation in three space dimensions ($d = 3$), imposing asymptotic boundary conditions and an initial condition 
\begin{subequations}
\label{eq:Equation1}
\begin{equation}
\label{eq:Equation1a}
\begin{gathered}
\begin{cases}
&\ii \, \partial_t \psi(x, t) = - \, \tfrac{1}{2} \, \Delta \, \psi(x, t) - \omega'(t) \, L(x) \, \psi(x, t) \\
&\qquad\qquad\qquad + \; V(x) \, \psi(x, t) + f\big(\absbig{\psi(x, t)}\big) \, \psi(x, t)\,, \\
&\psi(x, t_0) = \psi_0(x)\,, \quad (x, t) \in \RR^d \times (t_0, T)\,;
\end{cases}
\end{gathered}
\end{equation}
in view of a numerical illustration for the creation of vortices, it is convenient to consider the reduced two-dimensional model, which is included for $d = 2$. 
In~\eqref{eq:Equation1a}, $\psi: \RR^d \times [t_0,T] \to \CC: (x, t) \to \psi(x, t)$ denotes the complex-valued macroscopic wave function and $\Delta = \partial_{x_1}^2 + \dots + \partial_{x_d}^2$ the Laplace operator with respect to cartesian coordinates $x = (x_1, \dots, x_d) \in \RR^d$.
In general, the angular momentum rotation term has the form $\omega'(t) \cdot L(x)$ with direction of rotation given by the time derivative of a vector-valued function $\omega: [t_0, T] \to \RR^3$ and angular momentum operator defined by the cross product of position and gradient $L(x) = - \, \ii \, (x \times \nabla) \in \RR^3$.
Commonly, it is assumed that the axis of rotation coincides with the $x_3$-axis;
thus, it suffices to consider a scalar time-dependent function~$\omega: [t_0, T] \to \RR$ and set 
\begin{equation}
L(x) = - \, \ii \, \big(x_1 \, \partial_{x_2} - x_2 \, \partial_{x_1}\big)\,, \quad x = (x_1, \dots, x_d) \in \RR^d\,.
\end{equation}
Further decisive quantitites in~\eqref{eq:Equation1a} are the external potential $V: \RR^d \to \RR$ describing the confining trap and the nonlinear function $f: \RR \to \RR$ characterising the particle interaction; 
with regard to numerical simulations, we focus on the consideration of a quadratic potential involving positive weights 
\begin{equation}
V(x) = \tfrac{1}{2} \sum_{\ell=1}^{d} \gamma_{\ell}^2 \, x_{\ell}^2\,, \quad \gamma_{\ell} > 0\,, \quad \ell \in \{1, \dots, d\}\,, \quad x = (x_1, \dots, x_d) \in \RR^d\,,
\end{equation}
and a cubic nonlinearity 
\begin{equation}
f\big(\abs{\psi(x,t)}\big) \, \psi(x, t) = \vartheta \, \abs{\psi(x,t)}^2 \, \psi(x, t)\,, \quad \vartheta \in \RR\,, \quad (x,t) \in \RR^d \times [t_0, T]\,,
\end{equation}
describing the strength of short-range two-body interactions in the condensate.

\paragraph{Generalisation}
Our approach and the proposed exponential time integration methods are directly applicable to the case of a real-valued space-time-dependent potential
$V: \RR^d \times [t_0, T] \to \RR$ and a nonlinear term of the form~$f(|\psi(x,t)|)$.
The formulation~\eqref{eq:Equation1} in particular complies with~\cite{BaoMarahrensTangZhang2013,ZengZhang2009} for the special case of constant 
angular velocity
\begin{equation}
\omega(t) = \Omega \, t\,, \quad t \in [t_0, T]\,, \quad \Omega \in \RR\,; 
\end{equation}
\end{subequations}
for simplicity, we do not include additional integral terms describing long-range dipole-dipole interactions.

\paragraph{Reformulation in rotating Lagrangian coordinates}
In view of the introduction of space and time discretisation methods for the rotational Gross--Pitaevskii equation~\eqref{eq:Equation1}, we apply a transformation to the rotating frame
\begin{subequations}
\label{eq:Equation2}
\begin{equation}
\begin{gathered}
d = 2: \quad 
R(t) = \begin{pmatrix} \cos\big(\omega(t)\big) & \sin\big(\omega(t)\big) \\ - \sin\big(\omega(t)\big) & \cos\big(\omega(t)\big) \end{pmatrix}\,, \\
d = 3: \quad 
R(t) = \begin{pmatrix} \cos\big(\omega(t)\big) & \sin\big(\omega(t)\big) & 0 \\ - \sin\big(\omega(t)\big) & \cos\big(\omega(t)\big) & 0 \\ 0 & 0 & 1 \end{pmatrix}\,, \\
x = R(t) \, \xi\,, \quad t \in [t_0, T]\,,
\end{gathered}
\end{equation}
and attain the non-autonomous nonlinear Schr\"{o}dinger equation 
\begin{equation}
\begin{gathered}
\begin{cases}
&\ii \, \partial_t \varphi(\xi, t) = - \, \tfrac{1}{2} \, \Delta \, \varphi(\xi, t) + V\big(R(t) \, \xi\big) \, \varphi(\xi, t) + f\big(\absbig{\varphi(\xi, t)}\big) \, \varphi(\xi, t)\,, \\
&\varphi(\xi, t_0) = \psi_0(\xi)\,, \quad (\xi, t) \in \RR^d \times (t_0, T)\,;
\end{cases}
\end{gathered}
\end{equation}
\end{subequations}
detailed calculations based on the corresponding classical Hamiltonian are given in Appendix~\ref{sec:AppendixA}.
In case of a space-time-dependent potential $V: \RR^d \times [t_0, T] \to \RR$, the transformed equation instead involves $V\big(R(t) \, \xi, t\big)$.
\section{Quasi-Magnus exponential integrators}
\label{sec:Section3}
$\,$ \\[-1mm]

\noindent 
In this section, we introduce CFQM exponential integrators~\cite{AlvermannFehske2011,BlanesCasasThalhammer2017} and modified CFQM exponential integrators~\cite{BaderBlanesKopylov2018} for non-autonomous linear Schr{\"o}dinger equations of the form~\eqref{eq:IntroductionEqu4}; 
as indicated in Section~\ref{sec:Section1}, we favour operator splitting and Fourier spectral methods for their numerical realisation.   
It suffices to specify the numerical approximation for the first time step of length $h > 0$. 

\paragraph{Stepwise generalisation}
We first review fundamentals on interpolating functions and associated solution representations based on the Magnus expansion, since these basic principles are practical in the design of optimised high-order schemes; 
in this context, we restrict ourselves to the consideration of a system of non-autonomous linear ordinary differential equations 
\begin{subequations}
\label{eq:ODELinear} 
\begin{equation}
\begin{cases}
\ii \, u'(t) = H(t) \, u(t) = \big(A + B(t)\big) \, u(t)\,, \quad t \in (t_0, t_0 + h)\,, \\ 
u(t_0) = u_{0}\,,
\end{cases}
\end{equation}
involving a constant square complex matrix $A \in \CC^{M \times M}$ and a matrix-valued time-dependent function $B: [t_0, t_0 + h] \to \CC^{M \times M}$.
In view of~\eqref{eq:IntroductionEqu4}, it is justified to assume that the values of~$B$ commute
\begin{equation}
\label{eq:ODELinearCommutator} 
\big[B(s), B(t)\big] = B(s) \, B(t) - B(t) \, B(s) = 0\,, \quad s, t \in [t_0, t_0 + h]\,, 
\end{equation}
\end{subequations}
whereas the matrix commutator $[A, B(t)]$ in general does not vanish; 
for instance, diagonal matrices fulfill condition~\eqref{eq:ODELinearCommutator}. 
The formal extension of CFQM exponential integrators for non-autonomous linear ordinary differential equations to non-autonomous linear evolution equations is straightforward;
however, the generalisation of a rigorous convergence analysis for differential equations defined by matrices to differential equations defined by unbounded operators is more involved, see~\cite{BlanesCasasGonzalezThalhammer2018}.
Additional considerations lead to modified CFQM exponential integrators for non-autonomous linear Schr{\"o}dinger equations.

\paragraph{Fundamentals}
With regard to the introduction of schemes up to order six, which we expect to be favourable in efficiency in connection with non-autonomous linear and nonlinear Schr{\"o}dinger equations, we focus on the sixth-order Gauss--Legendre quadrature rule with nodes
\begin{subequations} 
\begin{equation}
\label{eq:GaussNodes}    
c_1 = \tfrac{1}{2} - \tfrac{\sqrt{15}}{10}\,, \quad c_2 = \tfrac{1}{2}\,, \quad c_3 = \tfrac{1}{2} + \tfrac{\sqrt{15}}{10}\,.
\end{equation}
The Lagrange polynomials through these nodes are defined by 
\begin{equation*}
\begin{gathered}
\nL_1: \RR \longrightarrow \RR: \sigma \longmapsto \frac{\sigma - c_2}{c_1 - c_2} \, \frac{\sigma - c_3}{c_1 - c_3}\,, \\
\nL_2: \RR \longrightarrow \RR: \sigma \longmapsto \frac{\sigma - c_1}{c_2 - c_1} \, \frac{\sigma - c_3}{c_2 - c_3}\,, \\
\nL_3: \RR \longrightarrow \RR: \sigma \longmapsto \frac{\sigma - c_1}{c_3 - c_1} \, \frac{\sigma - c_2}{c_3 - c_2}\,;
\end{gathered}
\end{equation*}
by construction, the following identities are fulfilled 
\begin{equation*}
\nL_k(c_k) = 1\,, \quad \nL_k(c_{\ell}) = 0\,, \quad k, \ell \in \{1, 2, 3\}\,, \quad \ell \neq k\,.
\end{equation*}
We consider the associated interpolating function 
\begin{equation*}
\begin{gathered}
\widehat{H}: \RR \longrightarrow \CC^{M \times M}:
t \longmapsto \nL_1\big(\tfrac{t - t_0}{h}\big) \, H_1 + \nL_2\big(\tfrac{t - t_0}{h}\big) \, H_2 + \nL_3\big(\tfrac{t - t_0}{h}\big) \, H_3\,, \\
H_k = H(t_0 + c_k h) = A + B_k\,, \quad B_k = B(t_0 + c_k h)\,, \quad k \in \{1, 2, 3\}\,,
\end{gathered}
\end{equation*}
which indeed satisfies the relation 
\begin{equation*}
\widehat{H}(t_0 + c_k h) = \nL_1(c_k) \, H_1 + \nL_2(c_k) \, H_2 + \nL_3(c_k) \, H_3 = H_k\,, \quad k \in \{1, 2, 3\}\,.
\end{equation*}
Substituting $\tau = t - t_0 - \frac{h}{2}$, a straightforward calculation yields the representation 
\begin{equation*}
\widehat{H}\big(\tau + t_0 + \tfrac{h}{2}\big)
= A + B_2 - \tfrac{\sqrt{15}}{3} \, \tfrac{\tau}{h} \, (B_1 - B_3) + \tfrac{10}{3} \, \tfrac{\tau^2}{h^2} \, (B_1 - 2 \, B_2 + B_3)\,;
\end{equation*}
consequently, Taylor series expansions
\begin{equation*}
B_3 - B_1 = \nO(h)\,, \quad B_3 - 2 \, B_2 + B_1 = \nO\big(h^2\big)\,, 
\end{equation*}
and the restriction to a suitably chosen interval verify the relations 
\begin{equation}
\label{eq:alpha123}
\begin{gathered}
\widehat{H}\big(\tau + t_0 + \tfrac{h}{2}\big)
= \ii \, \tfrac{1}{h} \, \Big(\alpha_1 + \alpha_2 \, \tfrac{\tau}{h} + \alpha_3 \, \tfrac{\tau^2}{h^2}\Big)\,, \quad 
\tau \in \big[- \tfrac{h}{2}, \tfrac{h}{2}\big]\,, \\
\alpha_1 = - \, \ii \, h \, \big(A + B_2\big) = \nO(h)\,, \\
\alpha_2 = - \, \ii \, h \, \tfrac{\sqrt{15}}{3} \, \big(B_3 - B_1\big) = \nO\big(h^2\big)\,, \\
\alpha_3 = - \, \ii \, h \, \tfrac{10}{3} \, \big(B_3 - 2 \, B_2 + B_1\big) = \nO\big(h^3\big)\,.
\end{gathered}
\end{equation}
For the following considerations, it is useful to employ a short notation for matrix commutators of the decisive quantities 
\begin{equation*}
[k \, \ell] = \big[\alpha_k, \alpha_{\ell}\big] = 
\begin{cases}
0\,, &k = \ell\,, \\
\nO\big(h^{k+\ell}\big)\,, \quad &k \neq \ell\,, 
\end{cases}
\quad k, \ell \in \{1, 2, 3\}\,;
\end{equation*}
with regard to the design of modified CFQM exponential integrators, we note that condition~\eqref{eq:ODELinearCommutator} implies $[2 \, 3] = 0$, whereas $\alpha_1$ in general does not commute with~$\alpha_2$ or~$\alpha_3$, respectively.
Accordingly, nested matrix commutators are denoted by 
\begin{equation*}
\begin{gathered}
[i \, j \ldots k \, \ell] = [\alpha_{i},[\alpha_{j}, [\ldots ,[\alpha_{k},\alpha_{\ell}]\ldots]]] = \nO\big(h^{i+j+\cdots+k+\ell}\big)\,, \\
i,j, \dots, k, \ell \in \{1, 2, 3\}\,. 
\end{gathered}
\end{equation*}
As shown in~\cite{Blanes2018}, the order of the underlying quadrature rule is inherited by the associated differential equation. 
More precisely, the solutions to the differential equation~\eqref{eq:ODELinear} and the related initial value problem 
\begin{equation*}
\begin{cases}
\ii \, v'(t) = \widehat{H}(t) \, v(t)\,, \quad t \in (t_0, t_0 + h)\,, \\ 
v(t_0) = u_{0}\,,
\end{cases}
\end{equation*}
coincide up to the desired order
\begin{equation*}
v(t_0 + h) - u(t_0 + h) = \nO\big(h^7\big)\,;
\end{equation*}
hence, it suffices to construct suitable approximations to~$v(t_0 + h)$.
For a non-autonomous linear differential equation, a formal solution representation by the matrix exponential based on the Magnus expansion~\cite{Magnus1954} is valid.
In the present situation, the first terms in this series expansion are given by 
\begin{equation*}
\begin{gathered}
v(t_0 + h) = \ee^{\Omega(h)} \, u_0\,, \\
\Omega(h)
= - \, \ii \int_{t_0}^{t_0 + h} \widehat{H}(t) \; \dd t
- \tfrac{1}{2} \int_{t_0}^{t_0 + h} \int_{t_0}^{t} \big[\widehat{H}(t), \widehat{H}(s)\big] \; \dd s \, \dd t + \dots\,;
\end{gathered}
\end{equation*}
due to the simple structure of the quadratic function~$\widehat{H}$, the arising integrals can be computed analytically and be expressed in terms of $\alpha_1, \alpha_2, \alpha_3$ as well as nested commutators, see also~\cite{BlanesCasasOteoRos2009,IserlesNorsett1999,MuntheKaasOwren1999}. 
As an illustration, we specify the leading contributions;
substituting $t = \tau + t_0 + \frac{h}{2}$ and inserting the representation~\eqref{eq:alpha123} yields
\begin{equation*}
- \, \ii \int_{t_0}^{t_0 + h} \widehat{H}(t) \; \dd t
= - \, \ii \int_{- \frac{h}{2}}^{\frac{h}{2}} \widehat{H}\big(\tau + t_0 + \tfrac{h}{2}\big) \; \dd \tau \\
= \alpha_1 + \tfrac{1}{12} \, \alpha_3\,, 
\end{equation*}
and a similar calculation verifies 
\begin{equation*}
- \, \tfrac{1}{2} \int_{t_0}^{t_0 + h} \int_{t_0}^{t} \big[\widehat{H}(t), \widehat{H}(s)\big] \; \dd s \, \dd t 
= - \, \tfrac{1}{12} \, [1 \, 2] + \tfrac{1}{240} \, [2 \, 3]\,.
\end{equation*}
Including terms up to order six
\begin{equation}
\label{eq:Omega6} 
\begin{split}
\Omega^{[6]} 
&= \alpha_{1} + \tfrac{1}{12} \, \alpha_{3} - \tfrac{1}{12} \, [1 \, 2] + \tfrac{1}{240} \, [2 \, 3] + \tfrac{1}{360} \, [1 \, 1 \, 3]
- \tfrac{1}{240} \, [2 \, 1 \, 2] \\
&\qquad + \tfrac{1}{720} \, [1 \, 1 \, 1 \, 2]\,, \\
\end{split}
\end{equation}
this procedure provides an appropriate approximation to~$v(t_0 + h)$ and hence to~$u(t_0 + h)$, that is 
\begin{equation}
\Omega^{[6]} - \Omega(h) = \nO\big(h^7\big)\,, \quad \ee^{\Omega^{[6]}} u_0 - u(t_0 + h) = \nO\big(h^7\big)\,;
\end{equation}
\end{subequations}
because of the symmetry of the Magnus expansion and the quadrature rule, only odd terms appear in $\Omega^{[6]}$.
The relations in~\eqref{eq:Omega6} are the starting point for the design of CFQM exponential integrators up to order six;
we recall that the additional condition $[2 \, 3] = 0$ is only used for modified CFQM exponential integrators.

\paragraph{CFQM exponential integrators} 
A $p$th-order CFQM exponential integrator for a non-autonomous linear ordinary differential equation of the form~\eqref{eq:ODELinear} is given by a composition of several matrix exponentials that comprise certain linear combinations of the values of the defining time-dependent matrix at specified nodes 
\begin{equation}
\label{eq:CFQMSchemeLinearODE} 
\begin{gathered}
u_1 = \exp\bigg(- \ii \, h \sum_{k=1}^{K} a_{Jk} \, H(t_0 + c_k h)\bigg) \, \times \cdots \\
\qquad\qquad \times \, \exp\bigg(- \ii \, h \sum_{k=1}^{K} a_{1k} \, H(t_0 + c_k h)\bigg) \, u_0\,, \\
u_1 - u(t_0 + h) = \nO\big(h^{p+1}\big)\,;  
\end{gathered}
\end{equation}
for instance, for fourth-order or sixth-order approximations based on three Gauss--Legendre quadrature nodes~\eqref{eq:GaussNodes} and positive integers $J \in \NN$, associated coefficients $a_{jk}$ for $(j,k) \in \{1, \dots, J\} \times \{1, 2, 3\}$ are determined such that the corresponding approximate solutions~$u_1$ agree with~$\ee^{\Omega^{[6]}} u_0$ up to terms of order four or six, respectively. 
Generally, these coefficients are obtained from a set of order conditions, which imply
\begin{equation}
\label{eq:MainCondition}
u_1 - \ee^{\Omega^{[p]}} u_0 = \nO\big(h^{p+1}\big)\,;  
\end{equation}
a practicable approach for the derivation of non-redundant conditions is described in~\cite{BlanesCasasThalhammer2017}.
We note that a higher number of exponentials $J \in \NN$ potentially offers the possibility to design optimised schemes.
The extension of~\eqref{eq:CFQMSchemeLinearODE} to non-autonomous linear Schr{\"o}dinger equations
\begin{equation*}
\begin{cases}
\ii \, u'(t) = H(t) \, u(t)\,, \quad t \in (t_0, t_0 + h)\,, \\ 
u(t_0) = u_{0}\,,
\end{cases}
\end{equation*}
and in particular to~\eqref{eq:IntroductionEqu4} relies on the composition of the solutions to the autonomous linear Schr{\"o}dinger equations 
\begin{subequations}
\label{eq:CFQMSchemeLinear} 
\begin{equation}
\begin{split}
&\begin{cases}
&\displaystyle \ii \, v_1'(t) = \sum_{k=1}^{K} a_{1k} \, H(t_0 + c_k h) \, v_1(t)\,, \quad t \in (t_0, t_0 + h)\,, \\
&v_1(t_0) = u_0\,, 
\end{cases} \\
&\qquad\qquad\quad\, \vdots \\ 
&\begin{cases}
&\displaystyle \ii \, v_J'(t) = \sum_{k=1}^{K} a_{Jk} \, H(t_0 + c_k h) \, v_J(t)\,, \quad t \in (t_0, t_0 + h)\,, \\
&v_J(t_0) = v_{J-1}(t_0 + h)\,, 
\end{cases} \\
\end{split}
\end{equation}
and yields a numerical approximation through
\begin{equation}
u_1 = v_J(t_0 + h)\,, \quad u_1 - u(t_0 + h) = \nO\big(h^{p+1}\big)\,.
\end{equation}
\end{subequations}
In our numerical experiments, we apply standard and optimised CFQM exponential integrators of orders $p = 2, 4, 6$;
for details, we refer to Section~\ref{sec:Section5}. 

\paragraph{Modified CFQM exponential integrators}
As indicated in Section~\ref{sec:Section1}, for non-autonomous linear Schr{\"o}dinger equations with Hamilton operator given by the Laplacian and a real-valued space-time-dependent potential, more efficient modified CFQM exponential integrators can be designed from a reduced set of order conditions;
more precisely, it is used that certain Lie-brackets of the $d$-dimensional Laplacian $\Delta = \partial_{\xi_1}^2 + \dots + \partial_{\xi_d}^2$ and a space-time-dependent multiplication operator $W: \RR^d \times [t_0, T] \to \RR$ simplify as follows 
\begin{equation}
\label{eq:LieBracket}
\begin{gathered}
\big[W(\cdot, t), W(\cdot, s)\big] = 0\,, \\
\big[\big[\Delta, W(\cdot, t) - W(\cdot, s)\big], W(\cdot, t) - W(\cdot, s)\big]
= 2 \, \sum_{\ell=1}^{d} \big(\partial_{\xi_{\ell}} \big(W(\cdot, t) - W(\cdot, s)\big)\Big)^2\,, \\
s, t \in [t_0, T]\,,
\end{gathered}
\end{equation}
see~\eqref{eq:IntroductionEqu4} and~\eqref{eq:IntroductionLieCommutator}. 
In order to sketch the approach, we reconsider the non-autonomous linear ordinary differential equation~\eqref{eq:ODELinear} and in particular assume that condition~\eqref{eq:ODELinearCommutator} holds;
a rigorous derivation in the context of partial differential equations will be the subject of future work. 
In the present situation, the matrix commutator~$[2 \, 3]$ vanishes;
as a consequence, relation~\eqref{eq:Omega6} and hence the consistency condition~\eqref{eq:MainCondition} reduce to 
\begin{equation}
\label{eq:Omega6Modified} 
\begin{gathered}
\Omega^{[6]} 
= \alpha_{1} + \tfrac{1}{1 \, 2} \, \alpha_{3} - \tfrac{1}{12} \, [1 \, 2] + \tfrac{1}{360} \, [1 \, 1 \, 3]
- \tfrac{1}{240} \, [2 \, 1 \, 2] + \tfrac{1}{720} \, [1 \, 1 \, 1 \, 2]\,, \\
u_1 - \ee^{\Omega^{[6]}} u_0 = \nO\big(h^7\big)\,.
\end{gathered}
\end{equation}
Moreover, due to the fact that~$\alpha_{2}$, $\alpha_{3}$, and the nested commutator~$[2 \, 1 \, 2]$ are of similar type, it is advantageous to incorporate terms of the form
\begin{equation*}
\exp\big(x_2 \, \alpha_{2} + x_3 \, \alpha_{3} + y \, [2 \, 1 \, 2]\big)
\end{equation*}
with appropriately chosen coefficients $x_2, x_3, y \in \RR$ into the composition~\eqref{eq:CFQMSchemeLinearODE}. 

\paragraph{Sixth-order scheme}
A particularly efficient sixth-order modified CFQM exponential integrator has been proposed in~\cite{BaderBlanesKopylov2018}; 
for a non-autonomous linear ordinary differential equation~\eqref{eq:ODELinear}, it is given by
\begin{equation*}
\begin{gathered}
u_1 
= \exp\big(- \, x_{12} \, \alpha_{2} + x_{13} \, \alpha_{3} + y \, [2 \, 1 \, 2]\big) \times \exp\big(x_{21} \, \alpha_{1} - x_{22} \, \alpha_{2} + x_{23} \, \alpha_{3}\big) \\
\qquad\qquad \times \exp\big(x_{21} \, \alpha_{1} + x_{22} \, \alpha_{2} + x_{23} \, \alpha_{3}\big) 
\times \exp\big(x_{12} \, \alpha_{2} + x_{13} \, \alpha_{3} + y \, [2 \, 1 \, 2]\big) \, u_0\,, \\
x_{12} = - \, x_{13} = - \, \tfrac{1}{60}\,, \quad x_{21} = \tfrac{1}{2}\,, \quad x_{22} = - \, \tfrac{2}{15}\,, \quad x_{23} = \tfrac{1}{40}\,, \quad 
y = \tfrac{1}{43200}\,.
\end{gathered}
\end{equation*}
Recalling~\eqref{eq:GaussNodes} and~\eqref{eq:alpha123}, a straightforward calculation shows that this composition of matrix exponentials is equivalent to 
\begin{equation*} 
u_1 = \ee^{- \, \ii \, h \, (\overline{B}_4 + h^2 \widetilde{B})} \, \ee^{- \, \ii \, \frac{h}{2} \, (A + \overline{B}_{3})} \, 
\ee^{- \, \ii \, \frac{h}{2} \, (A + \overline{B}_{2})} \, \ee^{- \, \ii \, h \, (\overline{B}_1 + h^2 \widetilde{B})} \, u_0\,;
\end{equation*}
the decisive quantities are defined by  
\begin{subequations}
\label{eq:ModifiedMagnus6Linear}
\begin{equation} 
\begin{gathered}
c_1 = \tfrac{1}{2} - \tfrac{\sqrt{15}}{10}\,, \quad c_2 = \tfrac{1}{2}\,, \quad c_3 = \tfrac{1}{2} + \tfrac{\sqrt{15}}{10}\,, \\ 
a_{11} = \tfrac{10 + \sqrt{15}}{180}\,, \quad a_{12} = - \, \tfrac{1}{9}\,, \quad a_{13} = \tfrac{10 - \sqrt{15}}{180}\,, \\
a_{21} = \tfrac{15 + 8\sqrt{15}}{90}\,, \quad a_{22} = \tfrac{2}{3}\,, \quad a_{23} = \tfrac{15 - 8 \sqrt{15}}{90}\,, \\
\overline{B}_1= a_{11} \, B(t_0 + c_1 h) + a_{12} \, B(t_0 + c_2 h) + a_{13} \, B(t_0 + c_3 h)\,, \\
\overline{B}_2= a_{21} \, B(t_0 + c_1 h) + a_{22} \, B(t_0 + c_2 h) + a_{23} \, B(t_0 + c_3 h)\,, \\
\overline{B}_3= a_{23} \, B(t_0 + c_1 h) + a_{22} \, B(t_0 + c_2 h) + a_{21} \, B(t_0 + c_3 h)\,, \\
\overline{B}_4= a_{13} \, B(t_0 + c_1 h) + a_{12} \, B(t_0 + c_2 h) + a_{11} \, B(t_0 + c_3 h)\,, \\
\widetilde{B} = \tfrac{1}{25920} \, \Big[\big[A, B(t_0 + c_3 h) - B(t_0 + c_1 h)\big], B(t_0 + c_3 h) - B(t_0 + c_1 h)\Big]\,.
\end{gathered}
\end{equation}
Accordingly to~\eqref{eq:CFQMSchemeLinear}, it is convenient to reformulate the extension to non-autonomous Schr{\"o}dinger equations~\eqref{eq:IntroductionEqu4} as 
\begin{equation}
\begin{split}
&\begin{cases}
&\ii \, v_1'(t) = \big(\overline{B}_1 + h^2 \, \widetilde{B}\big) \, v_1(t)\,, \quad t \in (t_0, t_0 + h)\,, \\
&v_1(t_0) = u_0\,, 
\end{cases} \\
&\begin{cases}
&\ii \, v_2'(t) = \tfrac{1}{2} \, \big(A + \overline{B}_{2}\big) \, v_2(t)\,, \quad t \in (t_0, t_0 + h)\,, \qquad\qquad\qquad \\
&v_2(t_0) = v_1(t_0+h)\,, 
\end{cases} \\
&\begin{cases}
&\ii \, v_3'(t) = \tfrac{1}{2} \, \big(A + \overline{B}_{3}\big) \, v_3(t)\,, \quad t \in (t_0, t_0 + h)\,, \\
&v_3(t_0) = v_2(t_0+h)\,, 
\end{cases} \\
&\begin{cases}
&\ii \, v_4'(t) = \big(\overline{B}_4 + h^2 \, \widetilde{B}\big) \, v_4(t)\,, \quad t \in (t_0, t_0 + h)\,, \\
&v_4(t_0) = v_3(t_0+h)\,;
\end{cases} 
\end{split}
\end{equation}
here, the matrix commutator is replaced by the Lie-bracket specified in~\eqref{eq:LieBracket}.
Numerical experiments indicate that this procedure yields a sixth-order approximation
\begin{equation}
u_1 = v_4(t_0 + h)\,, \quad u_1 - u(t_0 + h) = \nO\big(h^7\big)\,;
\end{equation}
\end{subequations}
a rigorous error analysis remains an open question.
\section{Application to rotational Gross-Pitaevskii equation}
$\,$ \\[-1mm]

\noindent
In this section, we describe our methodology for the space and time discretisation of the Gross--Pitaevskii equation with rotation term, see~\eqref{eq:Equation1} and~\eqref{eq:Equation2}. 
We indicate how to adapt the basic principles introduced in Section~\ref{sec:Section3} for non-autonomous linear Schr{\"o}dinger equations~\eqref{eq:IntroductionEqu4} to the more complex nonlinear case;
in particular, we generalise CFQM exponential integrators~\eqref{eq:CFQMSchemeLinear} as well as the sixth-order modified CFQM exponential integrator~\eqref{eq:ModifiedMagnus6Linear}.
We sketch the implementation of CFQM exponential integrators based on operator splitting methods, see also Section~\ref{sec:Section1};
the realisation of the sixth-order modified CFQM exponential integrator is in the same line. 
A rigorous justification and convergence analysis of CFQM exponential integrators for non-autonomous nonlinear Schr{\"o}dinger equations is out of the scope of the present contribution and will be the subject of future work; 
however, numerical illustrations presented in Section~\ref{sec:Section5} confirm the appropriateness of our approach.

\paragraph{Abstract evolution equation}
In order to identify similarities to non-autonomous linear ordinary differential equations~\eqref{eq:ODELinear} and non-autonomous linear Schr{\"o}dinger equations~\eqref{eq:IntroductionEqu4}, it is useful to consider the rotational Gross--Pitaevskii equation~\eqref{eq:Equation2} on a subinterval of length $h > 0$ and to reformulate it as abstract evolution equation
\begin{subequations} 
\label{eq:GPEAbstract}
\begin{equation} 
\label{eq:GPEAbstract1}
\begin{cases}
\ii \, u'(t) = H\big(t, u(t)\big) = \Big(A + B\big(t, u(t)\big)\Big) \, u(t)\,, \quad t \in (t_0, t_0 + h)\,, \\
u(t_0) = u_0\,;
\end{cases} 
\end{equation}
essentially, the values of the solution and the defining operators correspond to 
\begin{equation}
\begin{gathered}
u(t) = \varphi(\cdot, t)\,, \quad u_0 = \psi_0\,, \quad A = - \, \tfrac{1}{2} \, \Delta\,, \\
W(\cdot, t) = V\big(R(t) \, (\cdot)\big)\,, \quad
B\big(t, u(t)\big) = W(\cdot, t) + f\big(\absbig{\varphi(\cdot, t)}\big)\,. 
\end{gathered}
\end{equation}
\end{subequations}
We recall that arguments detailed in Section~\ref{sec:Section1} justify the consideration of a time-independent linear part~$A$ and a time-dependent nonlinear part~$B$; 
a decomposition into three operators as in~\eqref{eq:IntroductionEqu5} is no longer needed, see paragraph \emph{Strategy for full discretisation}~(ii).

\paragraph{Approach}
The calculus of Lie-derivatives provides powerful tools for the formal extension of the linear to the nonlinear case;
however, its application in connection with unbounded operators is not straightforward. 
Thus, as a first step, it is helpful to consider a non-autonomous nonlinear ordinary differential equation of the form~\eqref{eq:GPEAbstract1} involving a constant square complex matrix $A \in \CC^{M \times M}$ and a matrix-valued time-dependent function $B: [t_0, t_0 + h] \times \CC^M \to \CC^{M \times M}$.
Moreover, it is important to observe that the space-time-dependent potential and the nonlinearity act as multiplication operators and commute;
that is, an analogous relation to~\eqref{eq:ODELinearCommutator} holds with Lie-brackets replacing matrix commutators 
\begin{equation}
\big[B\big(s, u(t)\big), B\big(\sigma, u(\tau)\big)\big] = 0\,, \quad s, t, \sigma, \tau \in [t_0, t_0 + h]\,. 
\end{equation}
Suitable adaptations of the arguments given in Section~\ref{sec:Section3} yield the analogues of the fundamental expansions~\eqref{eq:Omega6} and~\eqref{eq:Omega6Modified}, and in succession, the analogoues of CFQM exponential integrators~\eqref{eq:CFQMSchemeLinear} and the modified CFQM exponential integrator~\eqref{eq:ModifiedMagnus6Linear}, which we describe next.

\paragraph{CFQM exponential integrators}
The generalisation of CFQM exponential integrators to the rotational Gross--Pitaevskii equation~\eqref{eq:GPEAbstract} requires the numerical solution of a sequence of autonomous Gross--Pitaevskii equations
\begin{subequations}
\label{eq:CFQMSchemeNonlinear}
\begin{equation}
\label{eq:CFQMSchemeNonlinear1}
\begin{cases}
&\displaystyle \ii \, v_j'(t) = \sum_{k=1}^{K} a_{jk} \, H\big(t_0 + c_k h, v_j(t)\big)\,, \quad t \in (t_0, t_0 + h)\,, \\
&v_j(t_0) = v_{j-1}(t_0 + h)\,,
\end{cases} \\
\end{equation}
where $v_1(t_0) = u_0$ and $j \in \{1, \dots, J\}$, yielding an approximation to the exact solution through
\begin{equation}
u_1 = v_J(t_0 + h)\,, \quad u_1 - u(t_0 + h) = \nO\big(h^{p+1}\big)\,,
\end{equation}
see also~\eqref{eq:CFQMSchemeLinearODE} and~\eqref{eq:CFQMSchemeLinear};
denoting 
\begin{equation*} 
b_j = \sum_{k=1}^{K} a_{jk}\,, \quad j \in \{1, \dots, J\}\,,
\end{equation*}
the nonlinear evolution equations arising in~\eqref{eq:CFQMSchemeNonlinear1} corresponds to
\begin{equation}
\label{eq:ResultingGPE}
\begin{cases}
&\displaystyle 
\ii \, \partial_t \varphi_j(\xi, t)
= - \, \tfrac{1}{2} \, b_j \, \Delta \, \varphi_j(\xi, t) + \sum_{k=1}^{K} a_{jk} \, W(\xi, t_0 + c_k h) \, \varphi_j(\xi, t) \\
&\qquad\qquad\qquad\quad + \; b_j \, f\big(\absbig{\varphi_j(\xi, t)}\big) \, \varphi_j(\xi, t)\,, \\
&\varphi_j(\xi, t_0) = \varphi_{j-1}(\xi, t_0 + h)\,, \quad (\xi, t) \in \RR^d \times (t_0, t_0 + h)\,, 
\end{cases}
\end{equation}
\end{subequations}
with $\varphi_j(\cdot, t) = v_j(t)$ for $t \in [t_0, t_0 + h]$ and $j \in \{1, \dots, J\}$.

\paragraph{Full discretisation}
Favourable space and time discretisation methods for~\eqref{eq:CFQMSchemeNonlinear1} and~\eqref{eq:ResultingGPE}, respectively, are provided by 
high-order operator splitting methods combined with Fourier spectral space discretisations, realised by fast Fourier techniques. 
They rely on the repeated numerical solution of subequations of the form 
\begin{equation*} 
\begin{split}
&(A)
\begin{cases}
\ii \, \tfrac{\dd}{\dd t} \, v_j^{(A)}(t) = b_j A \, v_j^{(A)}(t)\,, \quad t \in (t_0, t_0 + h)\,, \\[1mm]
v_j^{(A)}(t_0) \text{ given}\,,
\end{cases} \\
&(A)
\begin{cases}
\ii \, \partial_t \varphi_j^{(A)}(\xi, t) = - \, \tfrac{1}{2} \, b_j \, \Delta \, \varphi_j^{(A)}(\xi, t)\,, \\[1mm]
\varphi_j^{(A)}(\xi, t_0) \text{ given}\,, \quad (\xi, t) \in \RR^d \times (t_0, t_0 + h)\,, 
\end{cases} \\
&(B)
\begin{cases}
\displaystyle \ii \, \tfrac{\dd}{\dd t} \, v_j^{(B)}(t) = \sum_{k=1}^{K} a_{jk} \, B\big(t_0 + c_k h, v_j^{(B)}(t)\big) \, v_j^{(B)}(t)\,, \\[1mm]
v_j^{(B)}(t_0) \text{ given}\,,
\end{cases} \\
&(B)
\begin{cases}
\displaystyle \ii \, \partial_t \varphi_j^{(B)}(\xi, t) = \bigg(\sum_{k=1}^{K} a_{jk} \, W(\xi, t_0 + c_k h)
+ b_j \, f\big(\absbig{\varphi_j^{(B)}(\xi, t)}\big)\bigg) \, \varphi_j^{(B)}(\xi, t)\,, \\[1mm]
\varphi_j^{(B)}(\xi, t_0) \text{ given}\,, \quad (\xi, t) \in \RR^d \times (t_0, t_0 + h)\,;
\end{cases} 
\end{split}
\end{equation*}
we include the formulations as abstract evolution equations and the associated partial differential equations.  
The computational cost for the resolution of subequation~(A) amounts to the application of a fast Fourier transform and its inverse;
a characteristic invariance property permits to reduce the nonlinear subequation~(B) to the linear subequation
\begin{equation*} 
\begin{split}
&(B)
\begin{cases}
\displaystyle \ii \, \tfrac{\dd}{\dd t} \, v_j^{(B)}(t) = \sum_{k=1}^{K} a_{jk} \, B\big(t_0 + c_k h, v_j^{(B)}(t_0)\big) \, v_j^{(B)}(t)\,, \\[1mm]
v_j^{(B)}(t_0) \text{ given}\,,
\end{cases} \\
&(B)
\begin{cases}
\displaystyle \ii \, \partial_t \varphi_j^{(B)}(\xi, t) = \bigg(\sum_{k=1}^{K} a_{jk} \, W(\xi, t_0 + c_k h)
+ b_j \, f\big(\absbig{\varphi_j^{(B)}(\xi, t_0)}\big)\bigg) \, \varphi_j^{(B)}(\xi, t)\,, \\[1mm]
\varphi_j^{(B)}(\xi, t_0) \text{ given}\,, \quad (\xi, t) \in \RR^d \times (t_0, t_0 + h)\,,
\end{cases} 
\end{split}
\end{equation*}
and to resolve it by pointwise multiplication. 
Altogether, a $p$th-order approximation to $u(t_0 + h)$ results from a composition involving the flows of the subequations and suitably chosen coefficients 
\begin{equation*} 
\begin{gathered} 
\Phi(h)(u_0) - u(t_0 + h) = \nO\big(h^{p+1}\big)\,, \\
\Phi(h) = \Phi_J(h) \circ \cdots \circ \Phi_1(h)\,, \\
\Phi_j(h) = \Phi_j^{(B)}(\beta_s h) \circ \Phi_j^{(A)}(\alpha_s h) \circ \cdots \circ \Phi_j^{(B)}(\beta_1 h) \circ \Phi_j^{(A)}(\alpha_1 h)\,, \\
v_j^{(A)}(t_0 + \alpha_{\ell} h) = \Phi_j^{(A)}(\alpha_{\ell} h)\big(v_j^{(A)}(t_0)\big)\,, \quad
v_j^{(B)}(t_0 + \beta_{\ell} h) = \Phi_j^{(B)}(\beta_{\ell} h)\big(v_j^{(B)}(t_0)\big)\,, \\
\alpha_{\ell}, \beta_{\ell} \in \RR\,, \quad \ell \in \{1, \dots, s\}\,, \quad j \in \{1, \dots J\}\,.
\end{gathered} 
\end{equation*}
In essence, our implementation of CFQM exponential integrators is based on this formulation, i.e., our algorithm calls subroutines~$\Phi_j^{(A)}$ and~$\Phi_j^{(B)}$ depending on the coefficients~$(a_{jk})_{k=1}^K$ and~$b_j$, which yield approximations to the solutions of the subequations with adjusted time increments~$\alpha_{\ell} h$ and~$\beta_{\ell} h$.
Additional information on time-splitting Fourier spectral space discretisation methods for autonomous Gross--Pitaevskii equations is found for instance  in~\cite{Thalhammer2012,ThalhammerAbhau2012}, see also Section~\ref{sec:Section5}. 
Due to the presence of a trapping potential, the values of the macroscopic wave function tend to zero;
hence, a suitable restriction of the space domain and the application of the Fourier or, alternatively, Sine spectral method is justified, see~\cite{ThalhammerAbhau2012}.

\paragraph{Sixth-order scheme}
Based on the above considerations, it is straighforward to extend the sixth-order scheme~\eqref{eq:ModifiedMagnus6Linear} to the rotational Gross--Pitaevskii equation.
We employ the abbreviations 
\begin{subequations}
\label{eq:ModifiedMagnus6Nonlinear}
\begin{equation} 
\begin{gathered}
c_1 = \tfrac{1}{2} - \tfrac{\sqrt{15}}{10}\,, \quad c_2 = \tfrac{1}{2}\,, \quad c_3 = \tfrac{1}{2} + \tfrac{\sqrt{15}}{10}\,, \\ 
a_{11} = \tfrac{10 + \sqrt{15}}{180}\,, \quad a_{12} = - \, \tfrac{1}{9}\,, \quad a_{13} = \tfrac{10 - \sqrt{15}}{180}\,, \\
a_{21} = \tfrac{15 + 8\sqrt{15}}{90}\,, \quad a_{22} = \tfrac{2}{3}\,, \quad a_{23} = \tfrac{15 - 8 \sqrt{15}}{90}\,, \\
\overline{W}_1(\xi) = a_{11} \, W(\xi, t_0 + c_1 h) + a_{12} \, W(\xi, t_0 + c_2 h) + a_{13} \, W(\xi, t_0 + c_3 h)\,, \\
\overline{W}_2(\xi) = a_{21} \, W(\xi, t_0 + c_1 h) + a_{22} \, W(\xi, t_0 + c_2 h) + a_{23} \, W(\xi, t_0 + c_3 h)\,, \\
\overline{W}_3(\xi) = a_{23} \, W(\xi, t_0 + c_1 h) + a_{22} \, W(\xi, t_0 + c_2 h) + a_{21} \, W(\xi, t_0 + c_3 h)\,, \\
\overline{W}_4(\xi) = a_{13} \, W(\xi, t_0 + c_1 h) + a_{12} \, W(\xi, t_0 + c_2 h) + a_{11} \, W(\xi, t_0 + c_3 h)\,, \\
\widetilde{W}(\xi)  = \tfrac{1}{25920} \sum_{\ell=1}^{d} \big(\partial_{\xi_{\ell}} W(\xi, t_0 + c_3 h) - \partial_{\xi_{\ell}} W(\xi, t_0 + c_1 h)\big)^2\,,
\end{gathered}
\end{equation}
see also~\eqref{eq:LieBracket};
for the special case of a quadratic potential, the arising spatial derivatives are given in Appendix~\ref{sec:AppendixB}.
As a consequence, we obtain 
\begin{equation}
\begin{split}
&\begin{cases}
&\ii \, \partial_t \varphi_1(\xi, t) = \big(\overline{W}_1(\xi) + h^2 \, \widetilde{W}(\xi)\big) \, \varphi_1(\xi, t)\,, \\
&\varphi_1(\xi, t_0) = \psi_0(\xi)\,, \quad (\xi, t) \in \RR^d \times (t_0, t_0 + h)\,, 
\end{cases} \\
&\begin{cases}
&\ii \, \partial_t \varphi_2(\xi, t)
= \tfrac{1}{2} \, \big(- \tfrac{1}{2} \, \Delta + \overline{W}_{2}(\xi) + f(\abs{\varphi_2(\xi, t_0)})\big) \, \varphi_2(\xi, t)\,, \\
&\varphi_2(\xi, t_0) = \varphi_1(\xi, t_0+h)\,, \quad (\xi, t) \in \RR^d \times (t_0, t_0 + h)\,,
\end{cases} \\
&\begin{cases}
&\ii \, \partial_t \varphi_3(\xi, t)
= \tfrac{1}{2} \, \big(- \tfrac{1}{2} \, \Delta + \overline{W}_{3}(\xi) + f(\abs{\varphi_3(\xi, t_0)})\big) \, \varphi_3(\xi, t)\,, \\
&\varphi_3(\xi, t_0) = \varphi_2(\xi, t_0+h)\,, \quad (\xi, t) \in \RR^d \times (t_0, t_0 + h)\,,
\end{cases} \\
&\begin{cases}
&\ii \, \partial_t \varphi_4(\xi, t) = \big(\overline{W}_4 + h^2 \, \widetilde{W}\big) \, \varphi_4(\xi, t)\,, \\
&\varphi_4(\xi, t_0) = \varphi_3(\xi, t_0+h)\,, \quad (\xi, t) \in \RR^d \times (t_0, t_0 + h)\,;
\end{cases}
\end{split}
\end{equation}
\end{subequations}
the solutions to the first and fourth subequations are computed by pointwise multiplication, and the second and third subequation again correspond to autonomous Gross--Pitaevskii equations. 

\paragraph{Future work}
We believe that it is worth to investigate our approach in the context of ground state computations~\cite{DanailaProtas2017,ZengZhang2009}; 
the use of exponential integrators involving complex coefficients with positive real parts guarantees stability of the approximations. 
\section{Numerical results}
\label{sec:Section5}
$\,$ \\[-1mm]

\begin{figure}[t!]
\includegraphics[width=8cm]{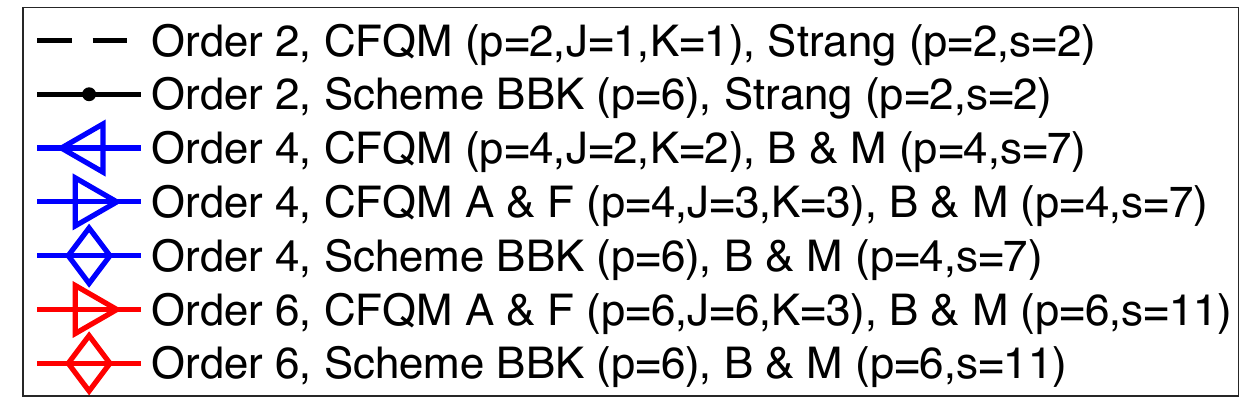} 
\caption{Overview of exponential time integration methods. Commutator-free quasi-Magnus exponential integrators (CFQM, order~$p$, number of compositions~$J$, number of quadrature nodes~$K$) or modified commutator-free quasi-Magnus exponential integrator (Scheme BBK, order $p = 6$), respectively, combined with operator splitting methods (order~$p$, number of compositions~$s$).}
\label{fig:FigureLegend}
\end{figure}

\noindent 
In this section, we illustrate the performance of the proposed full discretisation methods for the time-dependent Gross--Pitaevskii equation with rotation term in two and three space dimensions.
We provide numerical comparisons for CFQM exponential integrators and a modified CFQM exponential integrator; 
the numerical realisation relies on operator splitting and Fourier spectral space discretisation methods. 
We in particular implement the following standard and optimised schemes, see also Figure~\ref{fig:FigureLegend}.%
\footnote{For the convenience of the reader, the coefficients of the employed schemes are included in a file, available at 
\begin{quote}
\url{http://techmath.uibk.ac.at/mecht/MyHomepage/Research/MethodCoefficientsCFQMSplitting.rtf}.
\end{quote}}
\begin{enumerate}[(1)]
\item 
\textsc{Order 2, CFQM ($p = 2$, $J = 1$, $K = 1$), Strang ($p = 2$, $s = 2$)}:
The standard CFQM exponential integrator of order two involving a single Gaussian quadrature node and the solution of a single autonomous Gross--Pitaevskii equation, combined with the second-order Strang splitting method involving two compositions. 
\item 
\textsc{Order 4, CFQM ($p = 4$, $J = 2$, $K = 2$), B \& M ($p = 4$, $s = 7$)}:
The standard CFQM exponential integrator of order four involving two Gaussian quadrature nodes and the solution of two autonomous Gross--Pitaevskii equations, combined with an optimised fourth-order Runge--Kutta--Nyström splitting method (RKN$_74$) involving seven compositions proposed by \textsc{Blanes, Moan}~\cite{BlanesMoan2002Splitting}.   
\item 
\textsc{Order 4, CFQM A \& F ($p = 4$, $J = 3$, $K = 3$), B \& M ($p = 4$, $s = 7$)}:
An optimised CFQM exponential integrator of order four involving three Gaussian quadrature nodes and the solution of three autonomous Gross--Pitaevskii equations proposed by \textsc{Alvermann, Fehske}~\cite{AlvermannFehske2011}, again combined with the optimised fourth-order splitting method RKN$_74$.  
\item 
\textsc{Order 6, CFQM A \& F ($p = 6$, $J = 6$, $K = 3$), B \& M ($p = 6$, $s = 11$)}:
An optimised CFQM exponential integrator of order six involving three Gaussian quadrature nodes and the solution of six autonomous Gross--Pitaevskii equations proposed by \textsc{Alvermann, Fehske}~\cite{AlvermannFehske2011}, combined with an optimised sixth-order 
Runge--Kutta--Nyström splitting method (RKN$_{11}6$) involving eleven compositions proposed by \textsc{Blanes, Moan}~\cite{BlanesMoan2002Splitting}.
\item 
The sixth-order modified CFQM exponential integrator~\eqref{eq:ModifiedMagnus6Nonlinear} combined with the Strang splitting method or the optimised splitting methods RKN$_74$, RKN$_{11}6$, respectively, leading to approximations of orders $p = 2, 4, 6$;
with regard to the authors of~\cite{BaderBlanesKopylov2018}, we refer to~\eqref{eq:ModifiedMagnus6Nonlinear} as BBK scheme. 
\end{enumerate}

\begin{figure}[t!]
\includegraphics[width=6.2cm]{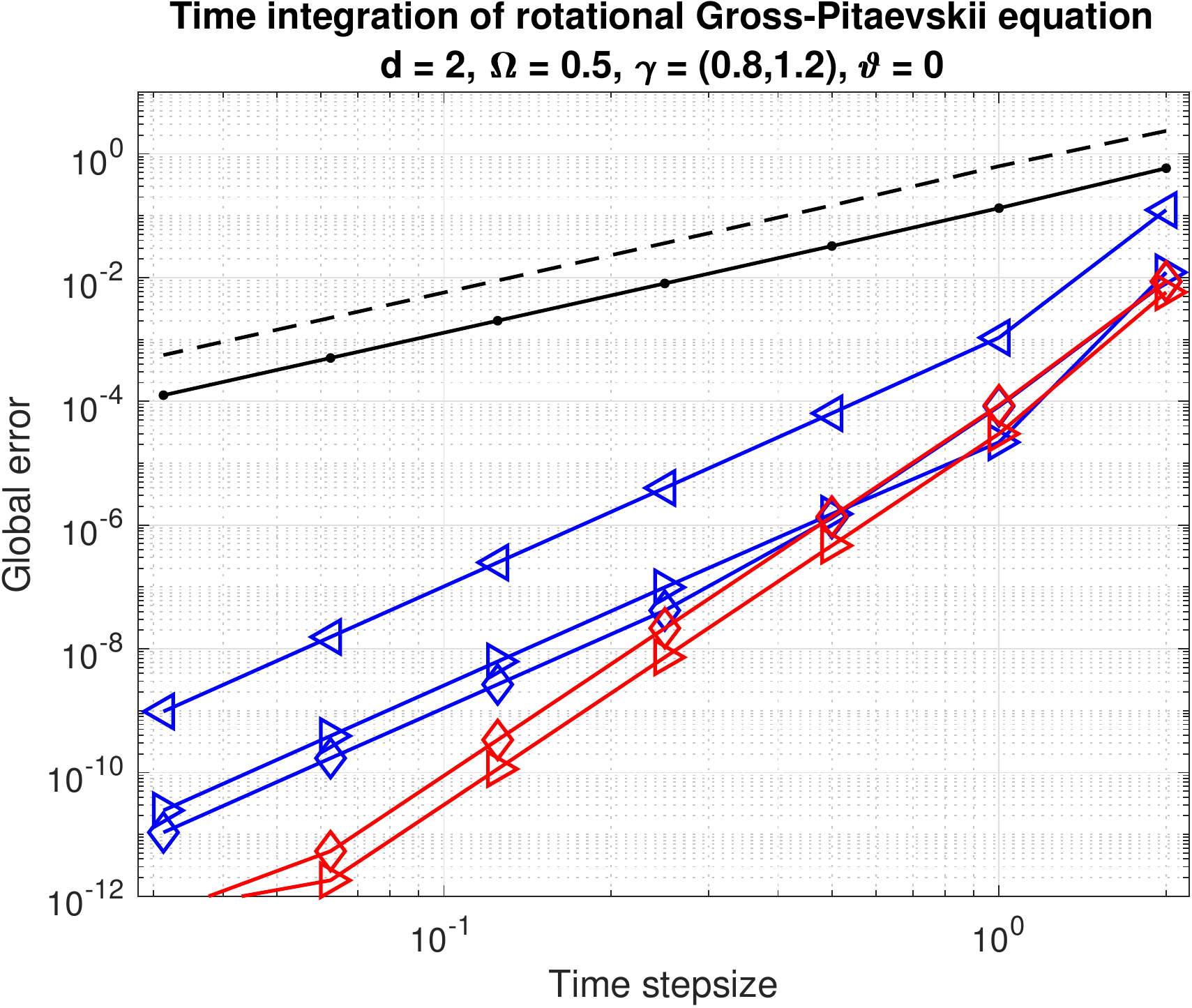} 
\includegraphics[width=6.2cm]{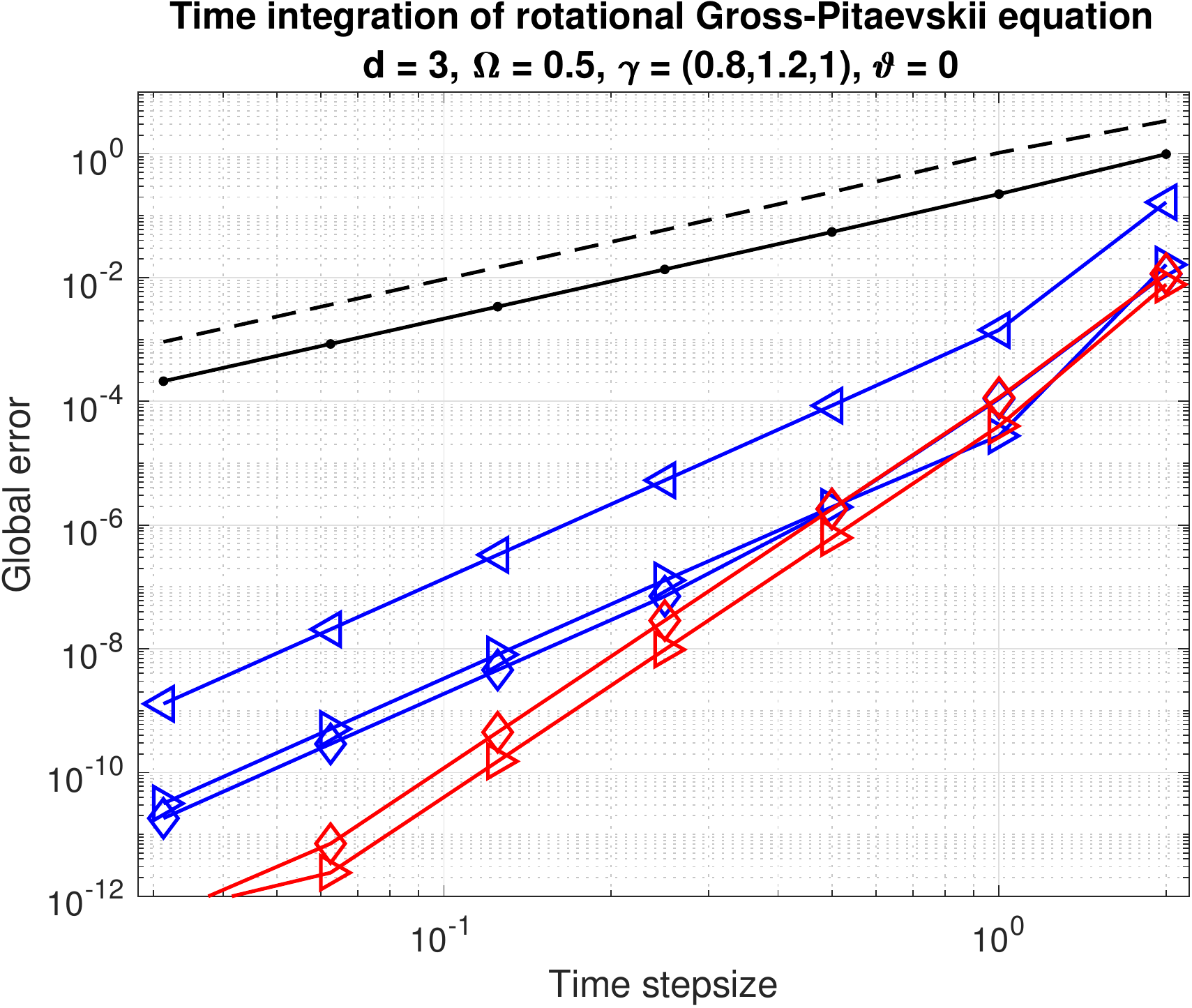} \\[2mm]
\includegraphics[width=6.2cm]{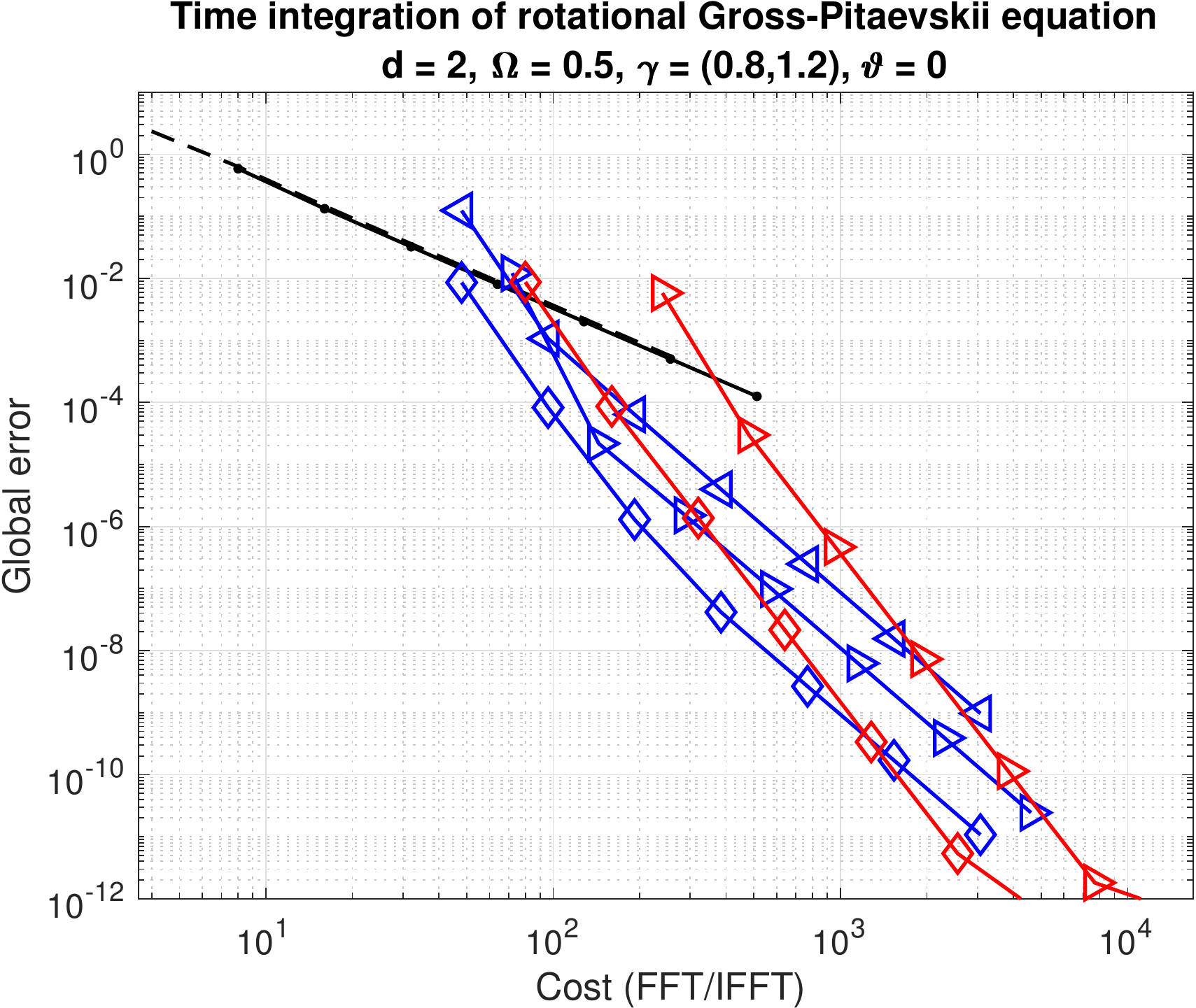} 
\includegraphics[width=6.2cm]{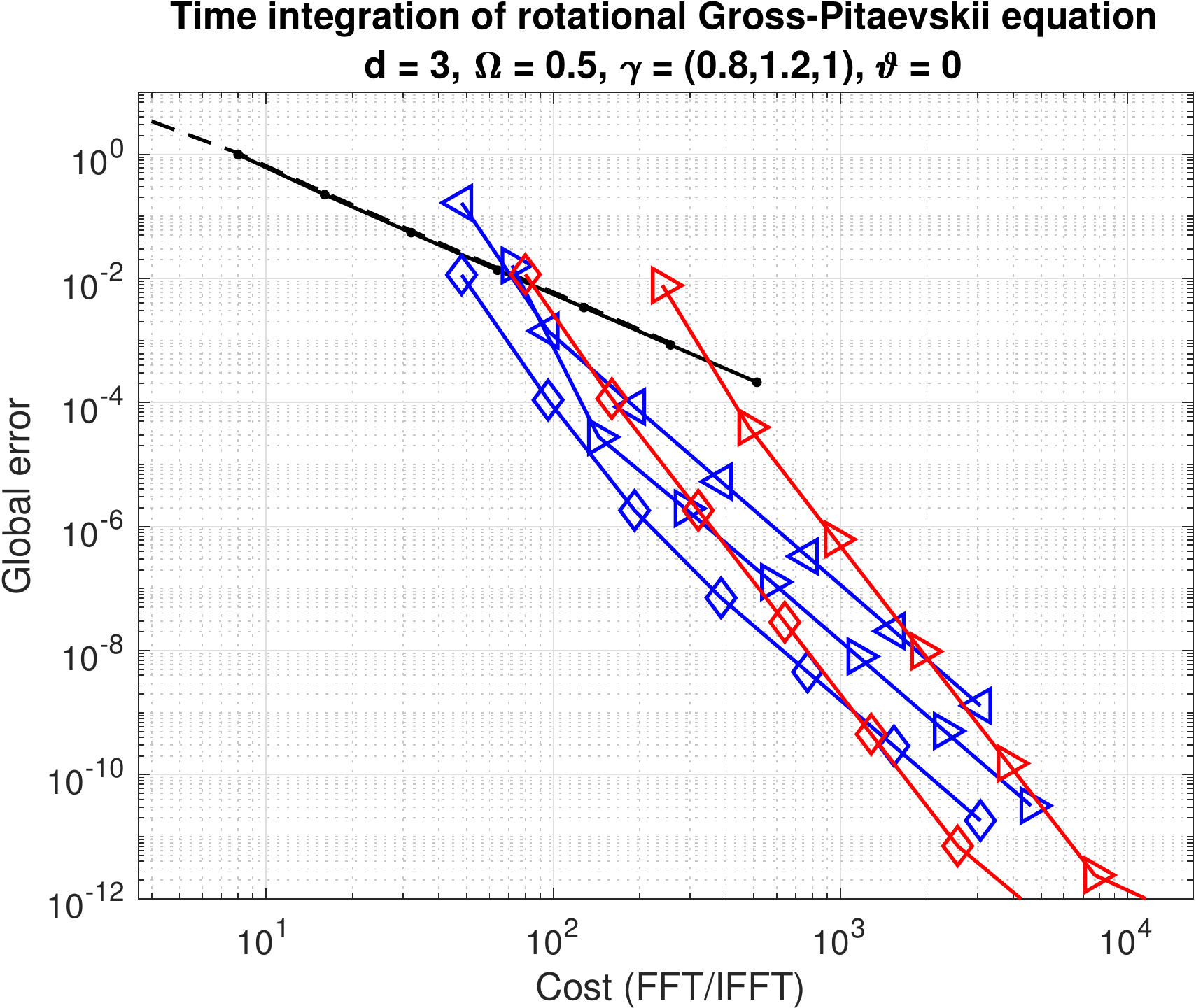} 
\caption{Time integration of non-autonomous linear test equation in two (left) and three (right) space dimensions by quasi-Magnus exponential integrators and operator splitting methods, see~\eqref{eq:Equation2}, \eqref{eq:TestEquation}, and Figure~\ref{fig:FigureLegend}. 
Global errors versus time stepsizes (first row) or number of fast Fourier transforms and inverse fast Fourier transforms (second row), respectively.}
\label{fig:FigureLinear}
\end{figure}

\begin{figure}[t!]
\includegraphics[width=6.2cm]{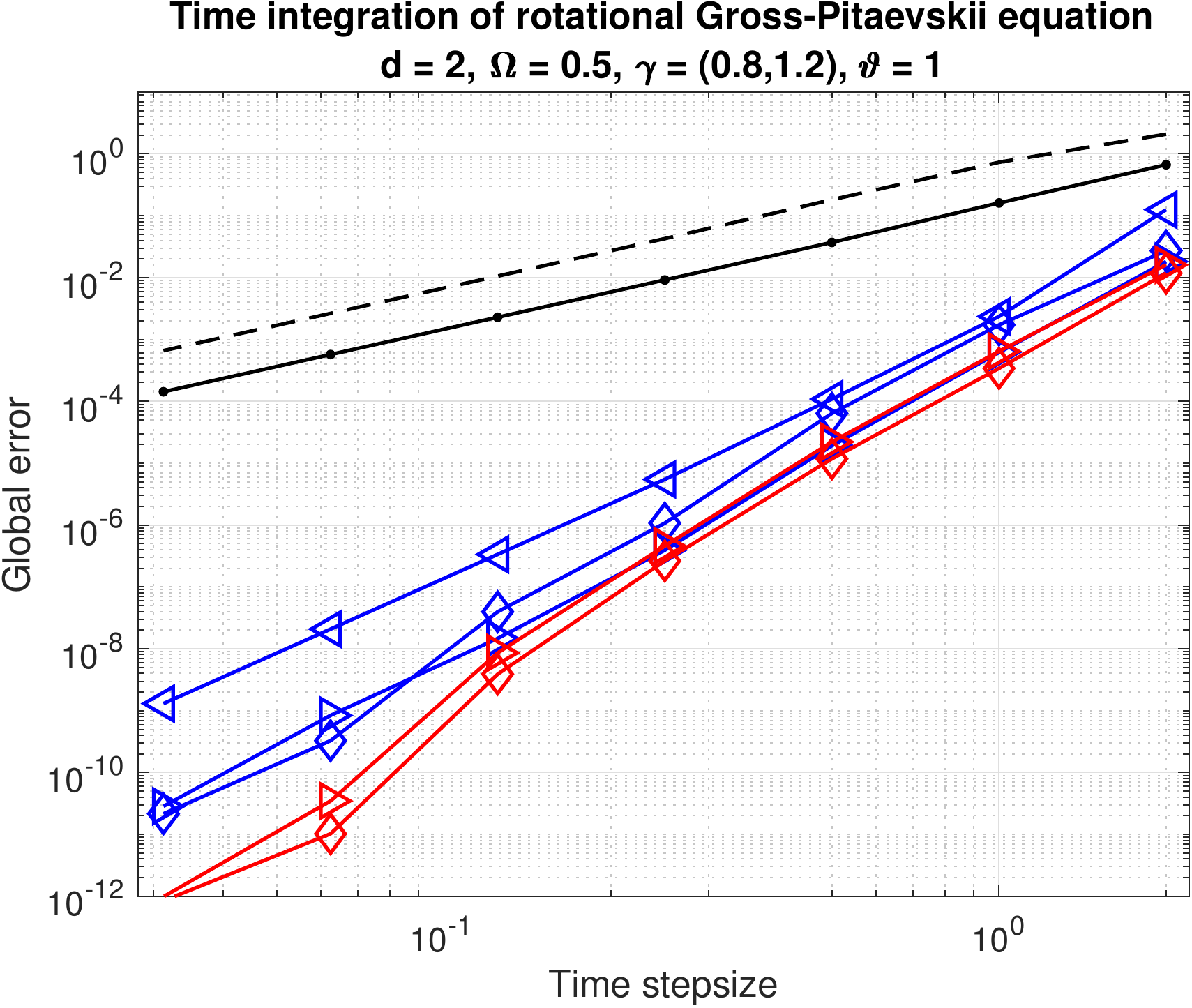} 
\includegraphics[width=6.2cm]{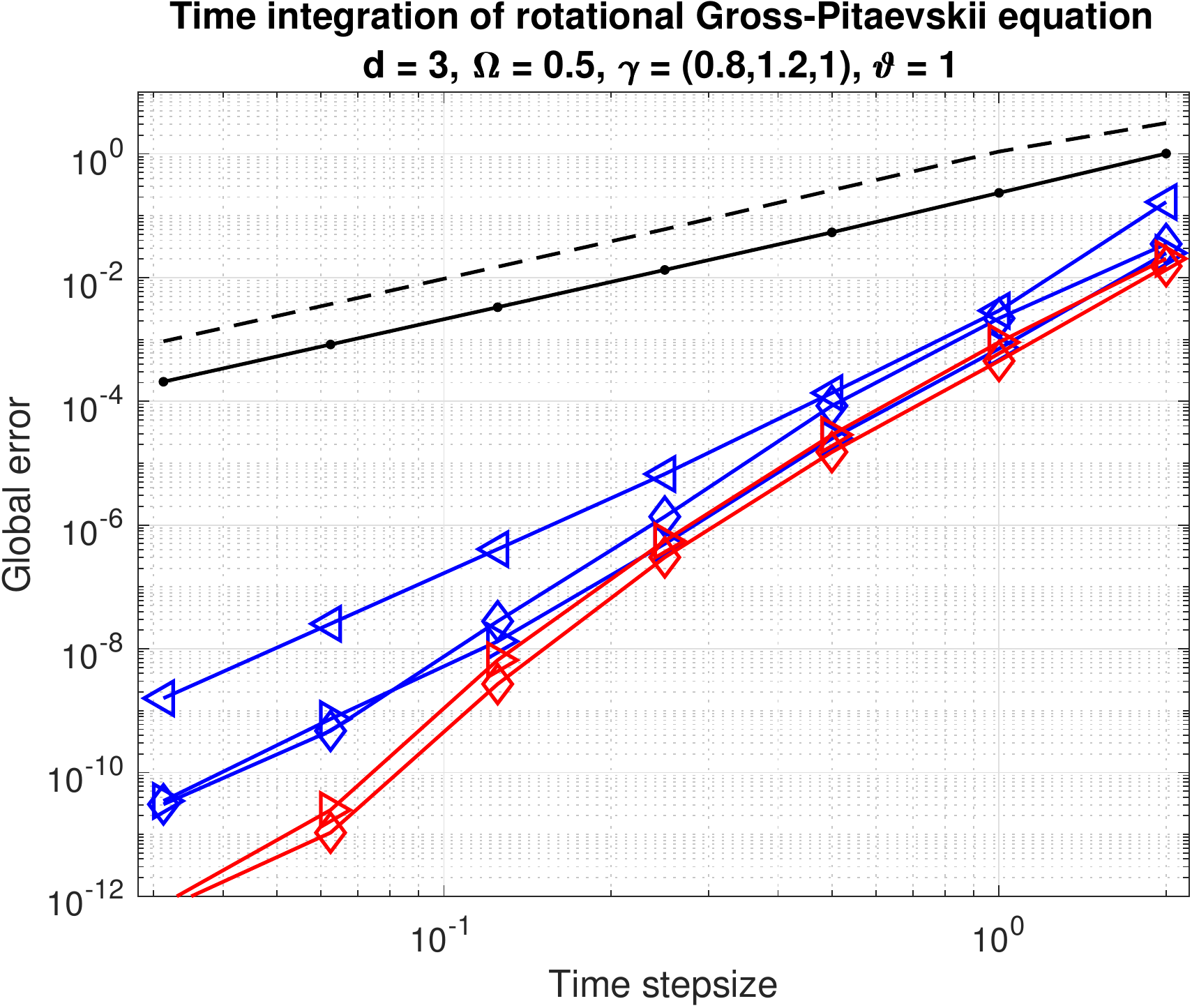} \\[2mm]
\includegraphics[width=6.2cm]{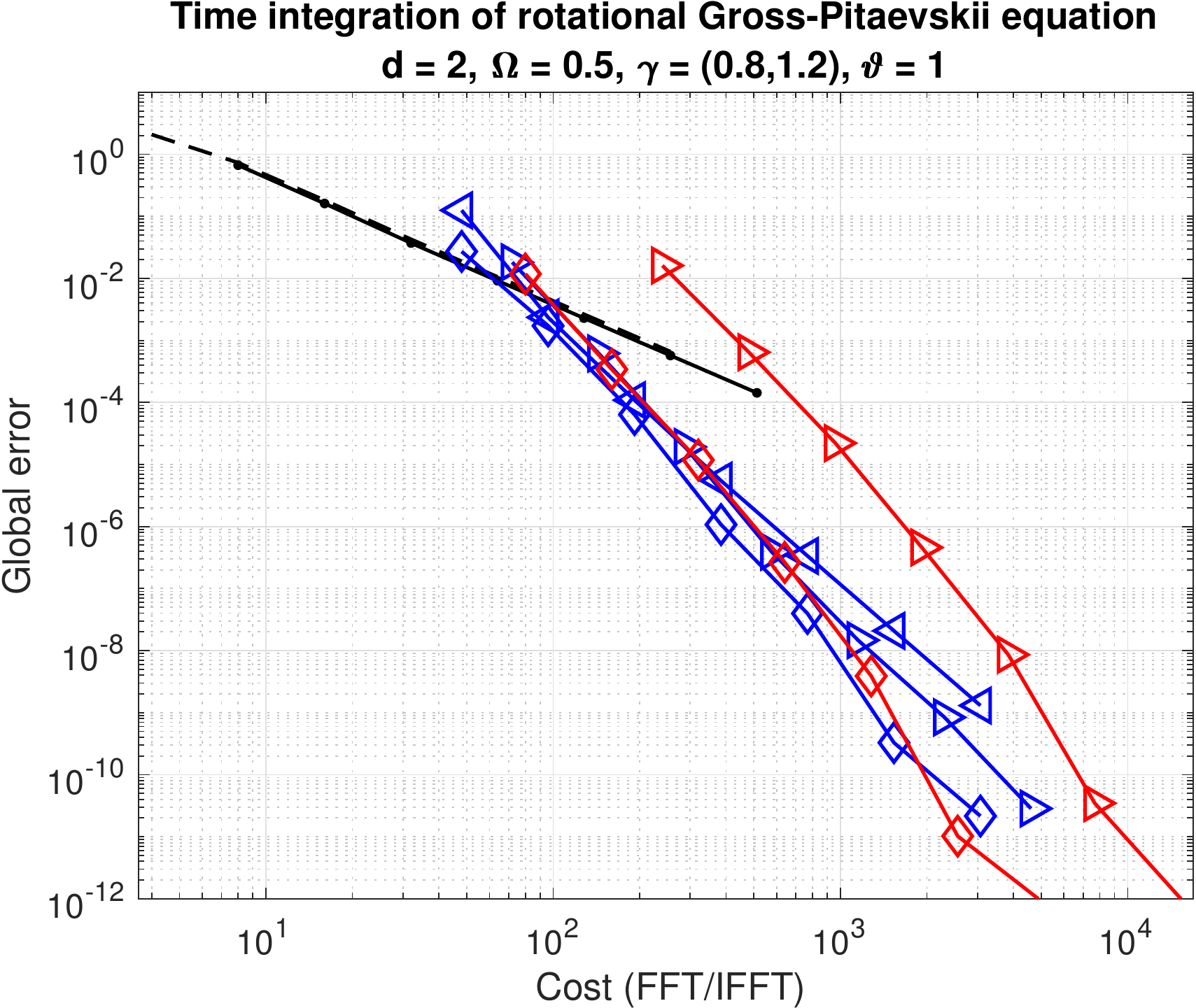}
\includegraphics[width=6.2cm]{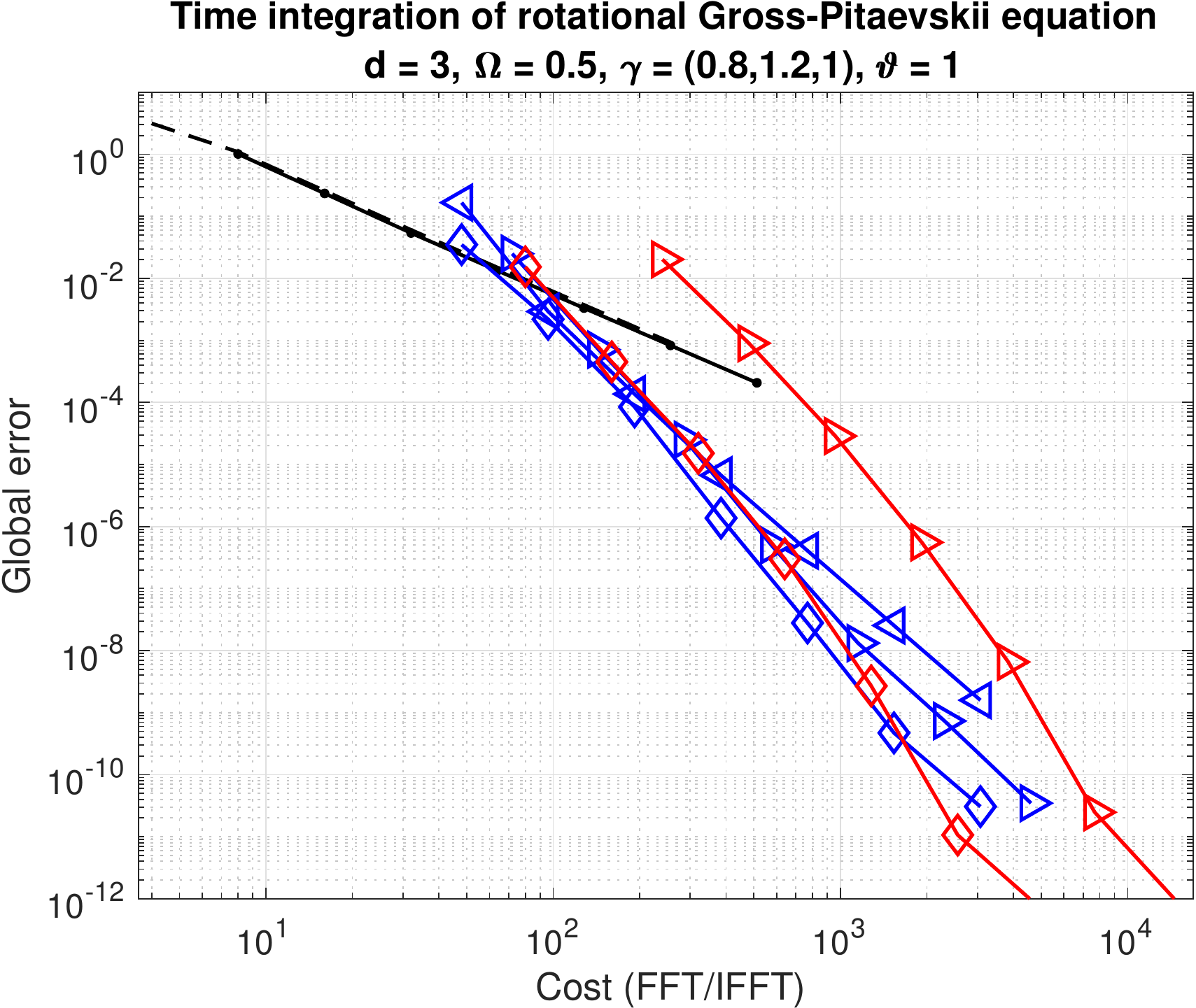} 
\caption{Time integration of non-autonomous nonlinear test equation in two (left) and three (right) space dimensions with $\vartheta = 1$ by quasi-Magnus exponential integrators and operator splitting methods, see~\eqref{eq:Equation2}, \eqref{eq:TestEquation}, and Figure~\ref{fig:FigureLegend}. 
Global errors versus time stepsizes (first row) or number of fast Fourier transforms and inverse fast Fourier transforms (second row), respectively.}
\label{fig:FigureNonlinear1}
\end{figure}

\begin{figure}[t!]
\includegraphics[width=6.2cm]{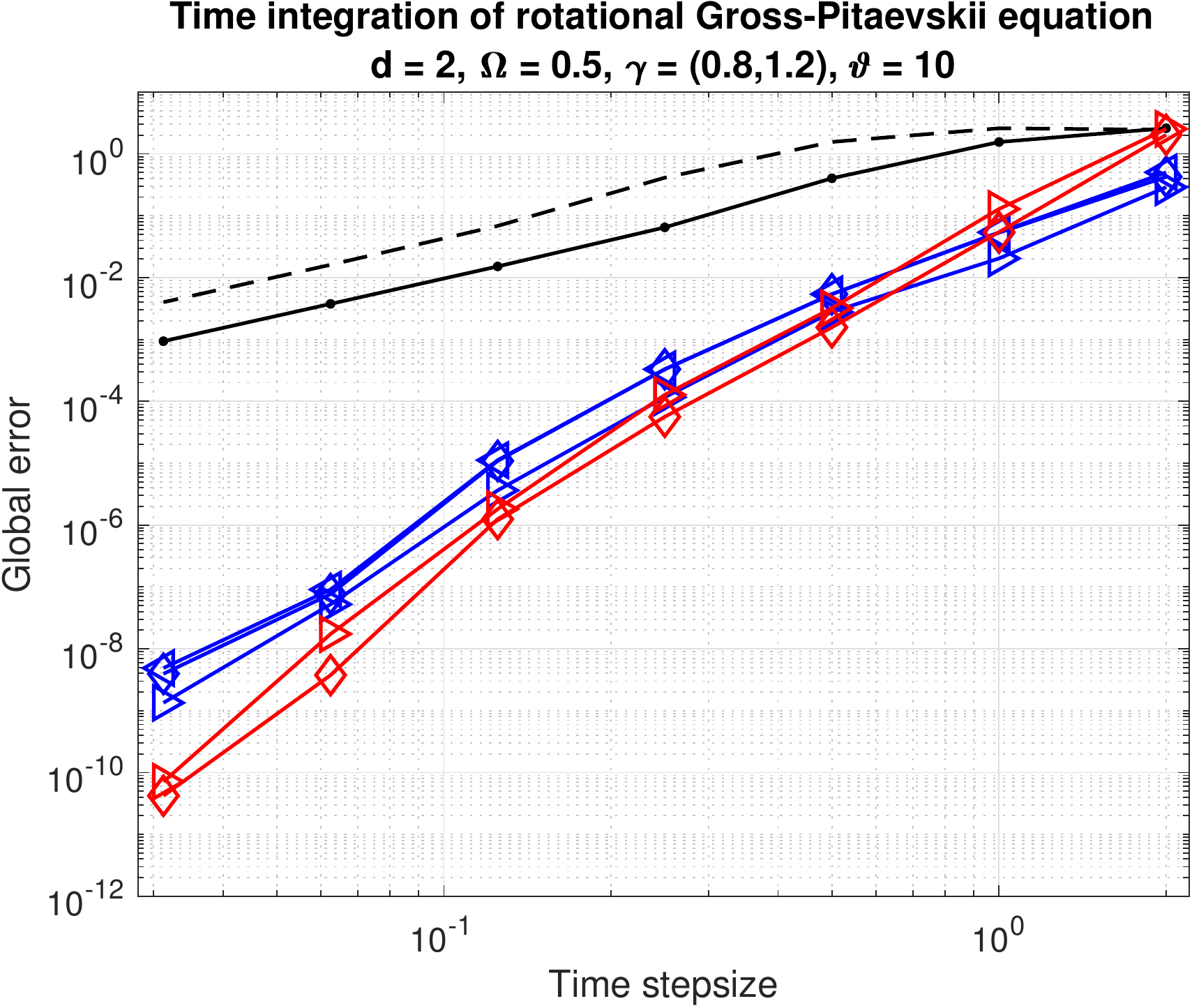} 
\includegraphics[width=6.2cm]{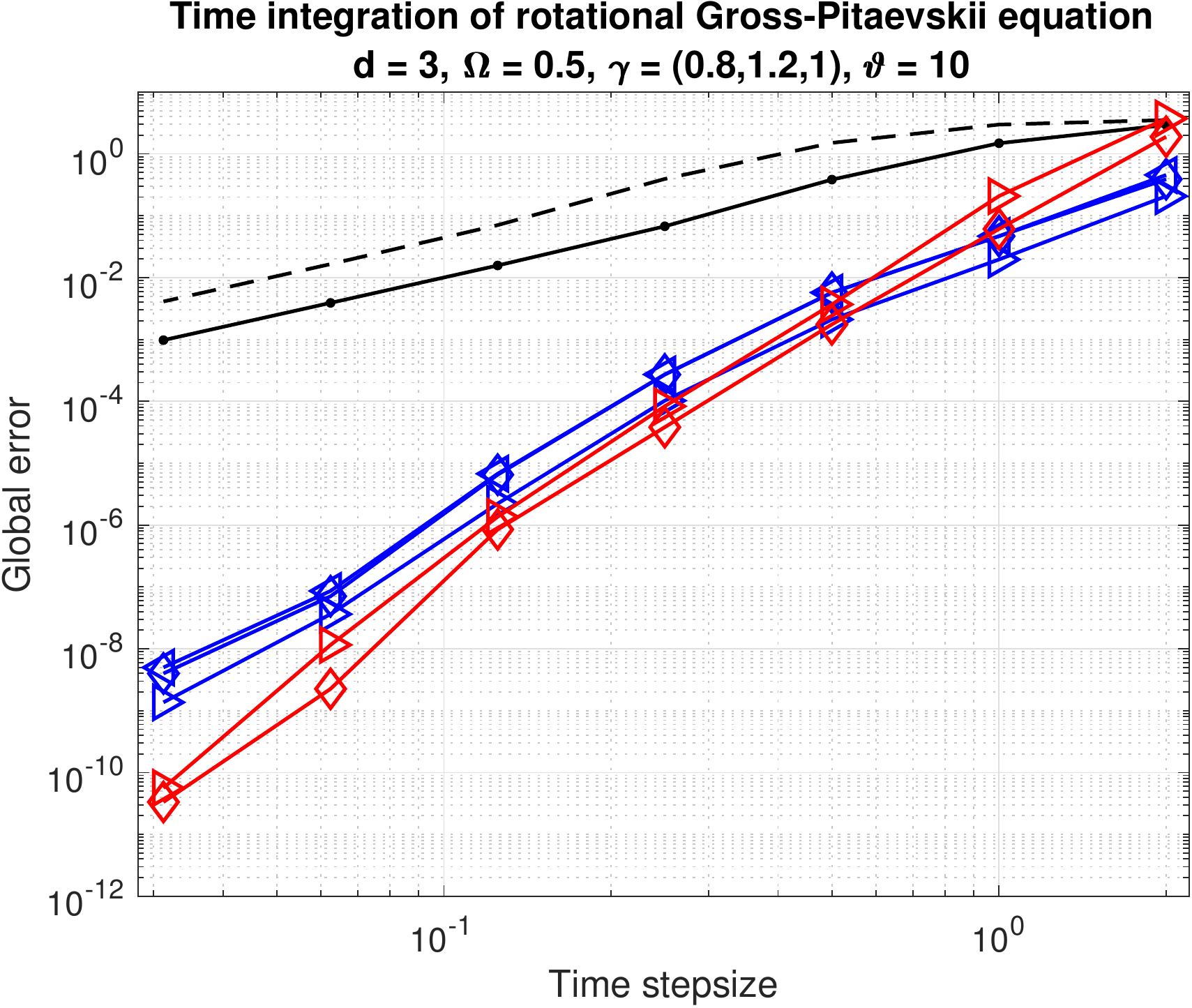} \\[2mm]
\includegraphics[width=6.2cm]{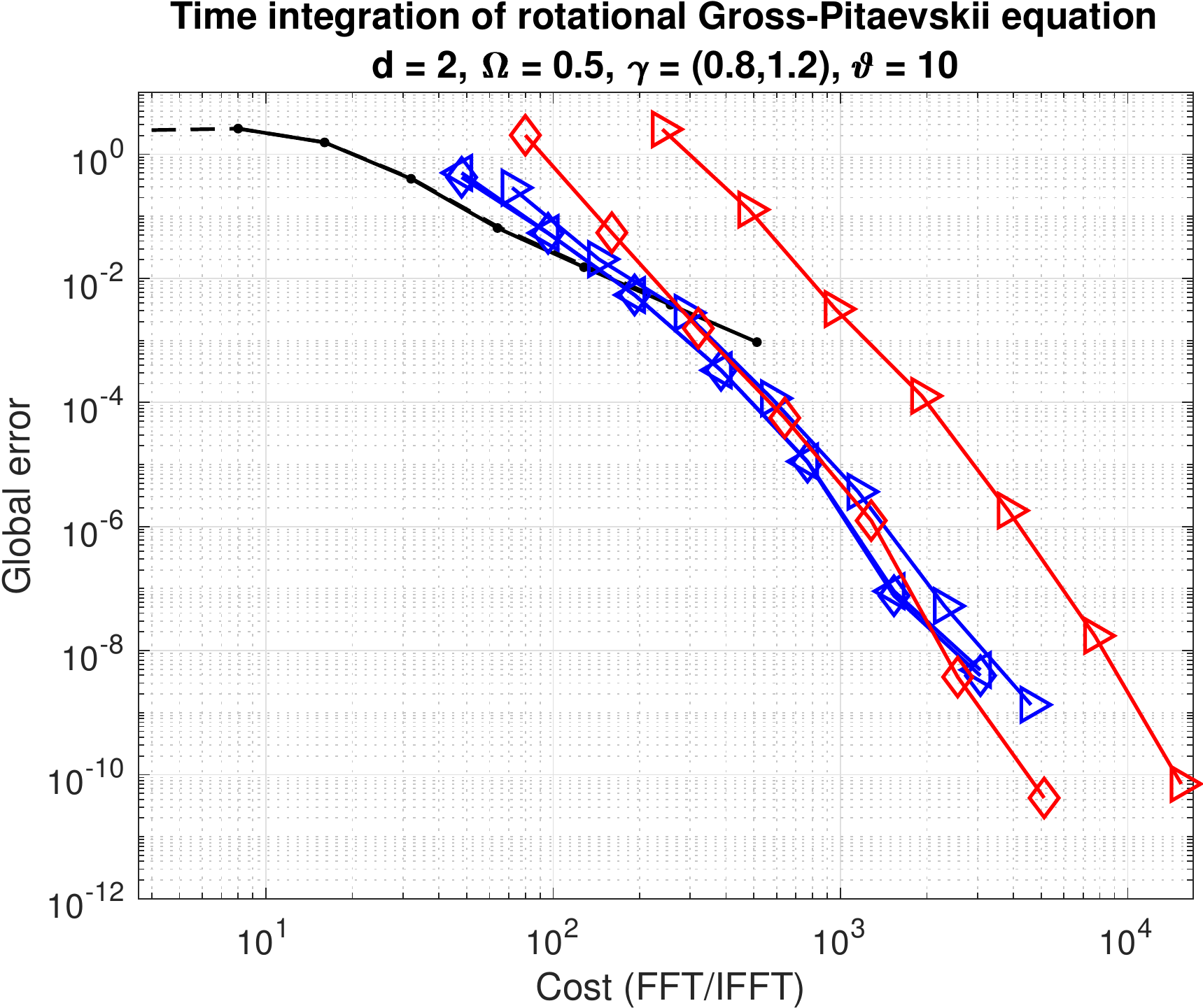}
\includegraphics[width=6.2cm]{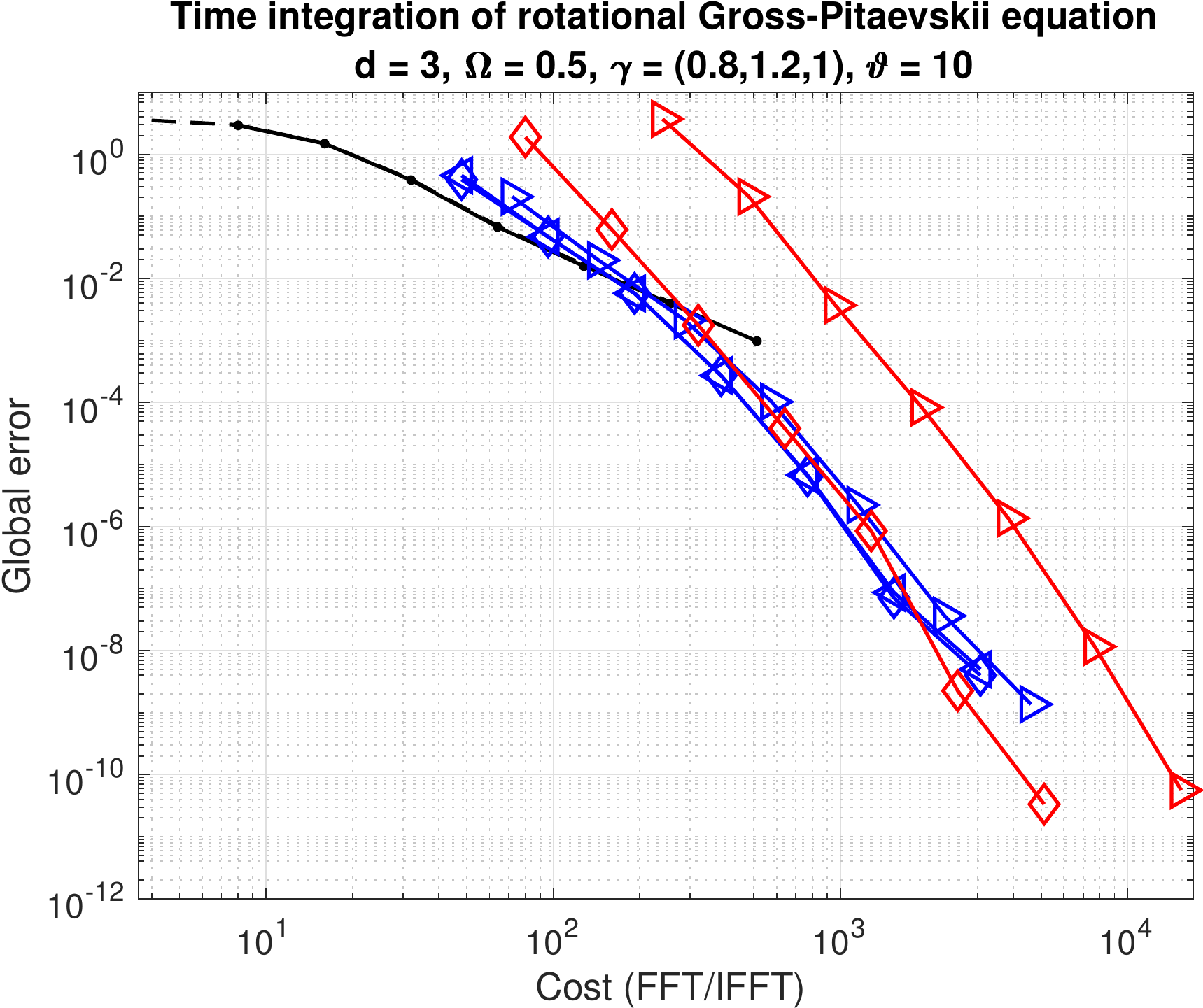}
\caption{Time integration of non-autonomous nonlinear test equation in two (left) and three (right) space dimensions with $\vartheta = 10$ by quasi-Magnus exponential integrators and operator splitting methods, see~\eqref{eq:Equation2}, \eqref{eq:TestEquation}, and Figure~\ref{fig:FigureLegend}. 
Global errors versus time stepsizes (first row) or number of fast Fourier transforms and inverse fast Fourier transforms (second row), respectively.}
\label{fig:FigureNonlinear2}
\end{figure}

\begin{figure}[ht!]
\includegraphics[width=8cm]{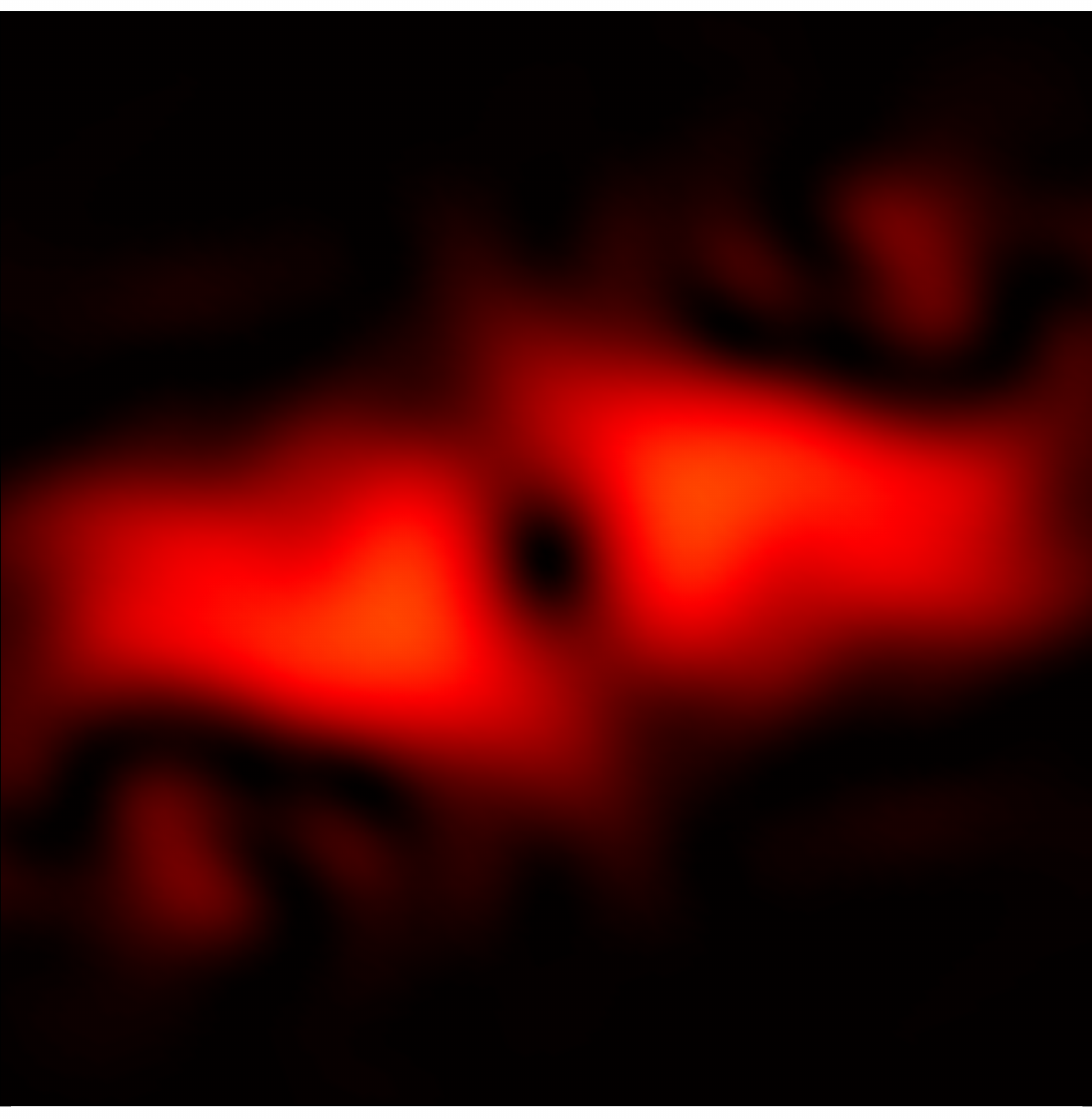} 
\caption{Time integration of two-dimensional rotational Gross--Pitaevskii equation~\eqref{eq:BEC} by sixth-order modified CFQM exponential integrator~\eqref{eq:ModifiedMagnus6Nonlinear} and optimised sixth-order splitting method proposed in~\cite{BlanesMoan2002Splitting}. 
Solution profile $\abs{\psi(x, t)}^2$ displayed for spatial section $x \in [- 5, 5]^2$ and time $t = 15$. 
Consistent results obtained by time-splitting generalised-Laguerre--Fourier spectral methods~\cite[Eq.\,(7)-(8), Fig.\,2]{HofstaetterKochThalhammer2014}.}
\label{fig:BEC}
\end{figure}

\paragraph{Linear and nonlinear test equations}
We consider the rotational Gross--Pitaevskii equation~\eqref{eq:Equation1}-\eqref{eq:Equation2} in two and three space dimensions, respectively, with initial and final times
\begin{subequations}
\label{eq:TestEquation}
\begin{equation}
t_0 = 0\,, \quad T = 4\,;
\end{equation}
in order to retain the nonstiff orders of convergence, we impose a smooth Gaussian-like initial state
\begin{equation}
\begin{gathered}
\psi_0(x) = \prod_{\ell=1}^{d} \ee^{- \frac{1}{2} \, \omega_{\ell}^2 \, x_{\ell}^2}\,, \quad x = (x_1, \dots, x_d) \in \RR^d\,, \\
d = 2: \quad \omega_1 = 1.1\,, \quad \omega_2 = 0.9\,, \\
d = 3: \quad \omega_1 = 1.1\,, \quad \omega_2 = 0.9\,, \quad \omega_3 = 1\,.
\end{gathered}
\end{equation}
The presence of a quadratic potential with weights 
\begin{equation}
\begin{gathered}
V(x) = \tfrac{1}{2} \sum_{\ell=1}^{d} \gamma_{\ell}^2 \, x_{\ell}^2\,, \quad x = (x_1, \dots, x_d) \in \RR^d\,, \\
d = 2: \quad \gamma_1 = 0.8\,, \quad \gamma_2 = 1.2\,, \\
d = 3: \quad \gamma_1 = 0.8\,, \quad \gamma_2 = 1.2\,, \quad \gamma_3 = 1\,,
\end{gathered}
\end{equation}
permits the restriction of the unbounded space domain to a cartesian product of bounded intervals and the use of a Fourier spectral space discretisation, realised by the fast Fourier transform and its inverse; 
we choose the following domains and numbers of equidistant space grid points 
\begin{equation}
\begin{gathered}
d = 2: \quad [- 10, 10] \times [- 10, 10]\,, \quad M = 64 \times 64\,, \\
d = 3: \quad [- 10, 10] \times [- 10, 10] \times [- 10, 10]\,, \quad M = 64 \times 64 \times 64\,.
\end{gathered}
\end{equation}
\end{subequations}
For the time integration of the rotational Gross--Pitaevskii equation, we apply the above listed exponential integrators of overall orders two, four, and six;
a numerical reference solution is determined by a sixth-order exponential integrator for the same spatial resolution and a refined time stepsize. 

In Figures~\ref{fig:FigureLinear}--\ref{fig:FigureNonlinear2}, we illustrate the efficiency of the considered exponential time integration methods for constant angular velocity $\Omega = 0.5$ and for increasing values of~$\vartheta$;
in the special case $\vartheta = 0$, the rotational Gross--Pitaevskii equation reduces to a non-autonomous linear Schr{\"o}dinger equation.
We in particular display the global errors versus the number of fast Fourier transforms and their inverse, which we consider to be the most costly part of the computations; 
the slopes of the lines reflect the temporal orders of convergence.
The obtained results for two and three space dimensions differ marginally.  
In all cases, it is observed that the modified CFQM exponential integrator~\eqref{eq:ModifiedMagnus6Nonlinear} is beneficial, in particular, when high accuracy is desired; 
this is explained by the fact that optimised CFQM integrators of orders four or six, respectively, require the solution of three or six, respectively, autonomous Gross--Pitaevskii equations in each time step, while the modified CFQM exponential integrator involves only two autonomous Gross--Pitaevskii equations. 

\paragraph{Solution profile}
In order to illustrate the dynamical behaviour of rotating Bose--Einstein condensates, we consider the rotational Gross--Pitaevskii equation~\eqref{eq:Equation1}-\eqref{eq:Equation2} with larger coupling constant and a particular initial state 
\begin{equation}
\label{eq:BEC}
\begin{gathered}
d = 2\,, \quad \Omega = 0.5\,, \quad \gamma_1 = 0.8\,, \quad \gamma_2 = 1.2\,, \quad \vartheta = 100\,, \quad t_0 = 0\,, \quad T = 15\,, \\
\psi_0(x) = \tfrac{1}{\sqrt{\pi}} \, (x_1 + \ii \, x_2) \, \ee^{- \frac{1}{2} \, (x_1^2 + x_2^2)}\,, \quad x = (x_1, x_2) \in \RR^2\,; 
\end{gathered}
\end{equation}
this problem has been studied in~\cite{HofstaetterKochThalhammer2014} in the context of the generalised-Laguerre--Fourier spectral method. 
As expected, the solution profile displayed in Figure~\ref{fig:BEC} is coincident with the result in~\cite[Fig.\,2]{HofstaetterKochThalhammer2014}. 
However, with respect to efficiency, our approach based on exponential integrators and Fourier spectral approximations outperforms the former approach based on the generalised-Laguerre--Fourier spectral method; 
moreover, its extension to the most relevant case of three space dimensions is straightforward. 
\section{Conclusions}
$\,$ \\[-1mm]

\noindent 
In this work, we proposed efficient exponential time integration methods for the Gross--Pitaevskii equation with rotation term in two and three space dimensions;
we used that a transformation to the rotating frame is straightforward by considering the associated classical Hamiltonian and employed a natural approach that combines CFQM exponential integrators with operator splitting and Fourier spectral space discretisation methods. 
On the basis of a sixth-order modified CFQM exponential integrator, we demonstrated that the special structure of the resulting nonlinear Schr{\"o}dinger equation allows the design of tailored schemes with reduced computational cost;
furthermore, in contrast to time-splitting generalised-Laguerre–Fourier–Hermite spectral methods, our approach permits the implementation of the physically relevant case of three space dimensions as a direct extension of the two-dimensional case.

Intrinsic advantages of CFQM exponential integrators and modified CFQM exponential integrators for non-autonomous linear differential equations are their favourable geometric properties, see~\cite{IserlesQuispel2018};
in the advanced context of non-autonomous linear Schr{\"o}dinger equations, this is reflected in an excellent stability behaviour in the underlying Lebesgue space of square-integrable functions.
In combination with optimised operator splitting methods and fast Fourier techniques, the proposed exponential integrators lead to efficient high-accuracy discretisations in space and time;
moreover, as described and illustrated in the present work, this methodology extends to non-autonomous nonlinear evolution equations.

Proceeding with a rigorous stability and error analysis, both for modified CFQM exponential integrators applied to non-autonomous linear Schr{\"o}dinger equations and for CFQM exponential integrators applied to non-autonomous nonlinear Schr{\"o}dinger equations, will be the subject of future work.
\section*{Acknowledgements}
$\,$ \\[-1mm]

\noindent 
Part of this work was developed during a research stay at the Wolfgang Pauli Institute Vienna; 
the authors are grateful to the director Norbert Mauser and the staff members for their support and hospitality. 
Philipp Bader, Sergio Blanes, and Fernando Casas acknowledge funding by the Ministerio de Economía y Competitividad (Spain) through project MTM2016-77660-P (AEI/FEDER, UE). 

\appendix
\section{Transformation to rotating Lagrangian coordinates}
\label{sec:AppendixA}
$\,$ \\[-1mm]

\paragraph{Simplification}
In this appendix, we include detailed calculations for the transformation of the time-dependent Gross--Pitaevskii equation with angular momentum rotation term~\eqref{eq:Equation1} by rotating Lagrangian coordinates~\eqref{eq:Equation2}.
For notational simplicity, we restrict ourselves to two space dimensions; 
setting $\xi_3 = x_3$, the extension to three space dimensions follows at once. 

\paragraph{Hamilton operator}
For the following derivation, it suffices to consider the Hamilton operator associated with the linear part in~\eqref{eq:Equation1}
\begin{equation*}
\begin{gathered}
\nH(x, t) = - \, \tfrac{1}{2} \, \Delta + \ii \, \omega'(t) \, \big(x_1 \, \partial_{x_2} - x_2 \, \partial_{x_1}\big) + V(x)\,, \\
x = (x_1, x_2) \in \RR^2\,, \quad t \in [t_0, T]\,.
\end{gathered}
\end{equation*}
Introducing position and momentum operators as well as compact vector and matrix notation, the Hamiltonian rewrites as  
\begin{equation}
\label{eq:nH}
\begin{gathered}
Q(x) = \big(x_1, x_2\big)^T\,, \quad P(x) = - \, \ii \, \big(\partial_{x_1}, \partial_{x_2}\big)^T\,, \quad 
J = \begin{pmatrix} 0 & 1 \\ - \, 1 & 0 \end{pmatrix} \in \RR^{2 \times 2}\,, \\
\nH(x, t) = \tfrac{1}{2} \, \big(P(x)\big)^T P(x) - \omega'(t) \, \big(Q(x)\big)^T J \, P(x) + V(x)\,, \\
x = (x_1, x_2) \in \RR^2\,, \quad t \in [t_0, T]\,.
\end{gathered}
\end{equation}

\paragraph{Classical Hamiltonian}
As common in classical mechanics, we denote the position and momentum of a particle by 
\begin{equation*}
q = \big(q_1, q_2\big)^T \in \RR^2\,, \quad p = \big(p_1, p_2\big)^T \in \RR^2\,;
\end{equation*}
based on the formal correspondence to the multiplication and differential operators in~\eqref{eq:nH}, the classical Hamiltonian
\begin{equation*}
\begin{gathered}
H(q, p, t) = \tfrac{1}{2} \, p^T p - \omega'(t) \, q^T J \, p + V(q)\,, \quad 
q, p \in \RR^2\,, \quad t \in [t_0, T]\,, 
\end{gathered}
\end{equation*}
is the analogue of the Hamilton operator associated with~\eqref{eq:Equation1}. 
Evidently, the partial derivatives are given by 
\begin{equation*}
\begin{gathered}
\partial_q H(q, p, t) = - \, \omega'(t) \, J \, p + V'(q)\,, \quad 
\partial_p H(q, p, t) = p - \omega'(t) \, J^T q = p + \omega'(t) \, J \, q\,, \\
q, p \in \RR^2\,, \quad t \in [t_0, T]\,, 
\end{gathered}
\end{equation*}
and the classical Hamiltonian system takes the form 
\begin{equation*}
\begin{split}
\begin{pmatrix} q'(t) \\ p'(t) \end{pmatrix}
&= \begin{pmatrix} \partial_p H\big(q(t), p(t), t\big) \\ - \, \partial_q H\big(q(t), p(t), t\big) \end{pmatrix} \\
&= \omega'(t) \begin{pmatrix} J & 0 \\ 0 & J \end{pmatrix} \begin{pmatrix} q(t) \\ p(t) \end{pmatrix} 
+ \begin{pmatrix} p(t) \\ - \, V'\big(q(t)\big) \end{pmatrix}\,, \quad t \in (t_0, T)\,;
\end{split}
\end{equation*}
employing the compact notation
\begin{equation*}
\begin{gathered}
z(t) = \begin{pmatrix} q(t) \\ p(t) \end{pmatrix} \in \RR^4\,, \quad \nJ = \begin{pmatrix} J & 0 \\ 0 & J \end{pmatrix} \in \RR^{4 \times 4}\, \\
A_{\text{rot}}(t) = \omega'(t) \, \nJ\,, \quad A_{\text{pot}}\big(z(t)\big) = \begin{pmatrix} p(t) \\- \, V'\big(q(t)\big)\end{pmatrix}\,, \quad t \in [t_0, T]\,,
\end{gathered}
\end{equation*}
we obtain the reformulation 
\begin{equation*}
z'(t) = A_{\text{rot}}(t) \, z(t) + A_{\text{pot}}\big(z(t)\big)\,, \quad t \in (t_0, T)\,.
\end{equation*}
Elementary calulations confirm the relations 
\begin{equation*}
\begin{gathered}
\ee^{\, \omega(t) J} 
= \begin{pmatrix} \cos\big(\omega(t)\big) & \sin\big(\omega(t)\big) \\[1mm] - \, \sin\big(\omega(t)\big) & \cos\big(\omega(t)\big) \end{pmatrix} = R(t)\,, \quad R'(t) = \omega'(t) \, J \, R(t)\,, \\
\tfrac{\dd}{\dd t} \, \ee^{\omega(t) \, \nJ} = \omega'(t) \, \nJ \, \ee^{\omega(t) \, \nJ} = A_{\text{rot}}(t) \, \ee^{\omega(t) \, \nJ}\,, \quad t \in (t_0, T)\,, 
\end{gathered}
\end{equation*}
which suggest the linear coordinate transform 
\begin{equation*}
\bar{z}(t) = \ee^{- \, \omega(t) \, \nJ} z(t)\,, \quad \bar{q}(t) = \ee^{- \, \omega(t) J} q(t)\,, \quad \bar{p}(t) = \ee^{- \, \omega(t) J} p(t)\,, \quad t \in [t_0, T]\,,
\end{equation*}
and lead to the transformed Hamiltonian system 
\begin{equation*}
\begin{gathered}
\bar{z}'(t) = \ee^{- \, \omega(t) \, \nJ} A_{\text{pot}}\big(\ee^{\omega(t) \, \nJ}\bar{z}(t)\big) 
= \begin{pmatrix} \bar{p}(t) \\ - \, \ee^{- \, \omega(t) J} \, V'\big(R(t) \, \bar{q}(t)\big) \end{pmatrix}\,, \quad t \in (t_0, T)\,;
\end{gathered}
\end{equation*}
observing that the associated Hamiltonian is given by 
\begin{equation*}
H(\bar{q}, \bar{p}, t) = \tfrac{1}{2} \, \bar{p}^T \bar{p} + V\big(R(t) \, \bar{q}\big)\,, \quad t \in [t_0, T]\,,
\end{equation*}
this implies the stated result~\eqref{eq:Equation2}.
\section{Auxiliary result}
\label{sec:AppendixB}
$\,$ \\[-1mm]

\noindent
Let $x = (x_1, x_2, x_3) \in \RR^3$, $\xi = (\xi_1, \xi_2, \xi_3) \in \RR^3$, and $t \in [t_0, T]$.  
In view of our numerical examples for the three-dimensional rotational Gross--Pitaevskii equation involving a quadratic potential
\begin{subequations}
\label{eq:QuadraticPotential}
\begin{equation}
V(x) = \tfrac{1}{2} \, \gamma_1^2 \, x_1^2 + \tfrac{1}{2} \, \gamma_2^2 \, x_2^2 + \tfrac{1}{2} \, \gamma_3^2 \, x_3^2\,, 
\end{equation}
we specify the arising space-time-dependent functions  
\begin{equation}
\begin{split}
R(t) \, \xi
&= \begin{pmatrix} \cos\big(\omega(t)\big) & \sin\big(\omega(t)\big) & 0 \\ - \sin\big(\omega(t)\big) & \cos\big(\omega(t)\big) & 0 \\ 0 & 0 & 1 \end{pmatrix} \!\!
\begin{pmatrix} \xi_1 \\ \xi_2 \\ \xi_3 \end{pmatrix} \\
&= \begin{pmatrix} \cos\big(\omega(t)\big) \, \xi_1 + \sin\big(\omega(t)\big) \, \xi_2 \\
- \sin\big(\omega(t)\big) \, \xi_1 + \cos\big(\omega(t)\big) \, \xi_2 \\ \xi_3 \end{pmatrix}, \\
W(\xi, t) &= V\big(R(t) \, \xi\big) \\
&= \tfrac{1}{2} \, \gamma_1^2 \, \Big(\cos\big(\omega(t)\big) \, \xi_1 + \sin\big(\omega(t)\big) \, \xi_2\Big)^2 \\
&\qquad + \tfrac{1}{2} \, \gamma_2^2 \, \Big(- \sin\big(\omega(t)\big) \, \xi_1 + \cos\big(\omega(t)\big) \, \xi_2\Big)^2
+ \tfrac{1}{2} \, \gamma_3^2 \, \xi_3^2\,, 
\end{split}
\end{equation}
see also~\eqref{eq:Equation1}-\eqref{eq:Equation2};
straightforward differentiation and simplification yield the relations
\begin{equation}
\begin{split}
\partial_{\xi_1} W(\xi, t)
&= \gamma_1^2 \, \cos\big(\omega(t)\big) \, \Big(\cos\big(\omega(t)\big) \, \xi_1 + \sin\big(\omega(t)\big) \, \xi_2\Big) \\
&\qquad - \gamma_2^2 \, \sin\big(\omega(t)\big) \, \Big(- \sin\big(\omega(t)\big) \, \xi_1 + \cos\big(\omega(t)\big) \, \xi_2\Big) \\
&= \gamma_1^2 \, \cos^2\big(\omega(t)\big) \, \xi_1 + \gamma_2^2 \, \sin^2\big(\omega(t)\big) \, \xi_1 \\
&\qquad+ \big(\gamma_1^2 - \gamma_2^2\big) \, \sin\big(\omega(t)\big) \, \cos\big(\omega(t)\big) \, \xi_2\,, \\
\partial_{\xi_2} W(\xi, t)
&= \gamma_1^2 \, \sin\big(\omega(t)\big) \, \Big(\cos\big(\omega(t)\big) \, \xi_1 + \sin\big(\omega(t)\big) \, \xi_2\Big) \\
&\qquad + \gamma_2^2 \, \cos\big(\omega(t)\big) \, \Big(- \sin\big(\omega(t)\big) \, \xi_1 + \cos\big(\omega(t)\big) \, \xi_2\Big) \\
&= \gamma_2^2 \, \cos^2\big(\omega(t)\big) \, \xi_2 + \gamma_1^2 \, \sin^2\big(\omega(t)\big) \, \xi_2 \\
&\qquad + \big(\gamma_1^2 - \gamma_2^2\big) \, \sin\big(\omega(t)\big) \, \cos\big(\omega(t)\big) \, \xi_1\,, \\
\partial_{\xi_3} W(\xi, t)
&= \gamma_3^2 \, \xi_3\,.
\end{split}
\end{equation}
\end{subequations}

\bigskip
\end{document}